\pdfoutput=1
\RequirePackage{ifpdf}
\ifpdf 
\documentclass[pdftex]{sigma}
\else
\documentclass{sigma}
\fi

\usepackage{here}	
\usepackage{myGlobalMacros}

\numberwithin{equation}{section}

\newtheorem{Theorem}{Theorem}[section]
\newtheorem{Corollary}[Theorem]{Corollary}
\newtheorem{Lemma}[Theorem]{Lemma}
\newtheorem{Proposition}[Theorem]{Proposition}
\newtheorem{Conjecture}[Theorem]{Conjecture}
\newtheorem{result}[Theorem]{Main results}
 { \theoremstyle{definition}
\newtheorem{Definition}[Theorem]{Definition}
\newtheorem{Example}[Theorem]{Example}
\newtheorem{Remark}[Theorem]{Remark} }

\newcommand{\abs}[1]{\lvert#1\rvert}

\begin{document}

\allowdisplaybreaks

\newcommand{\arXivNumber}{1301.7632}

\renewcommand{\PaperNumber}{067}

\FirstPageHeading

\ShortArticleName{Minuscule Schubert Varieties and Mirror Symmetry}

\ArticleName{Minuscule Schubert Varieties and Mirror Symmetry}

\Author{Makoto MIURA}

\AuthorNameForHeading{M.~Miura}

\Address{Korea Institute for Advanced Study,\\ 85 Hoegiro, Dongdaemun-gu, Seoul, 130-722, Republic of Korea}
\Email{\href{mailto:miura@kias.re.kr}{miura@kias.re.kr}}

\ArticleDates{Received August 23, 2016, in f\/inal form August 16, 2017; Published online August 23, 2017}

\Abstract{We consider smooth complete intersection Calabi--Yau 3-folds in minuscule Schubert varieties, and study their mirror symmetry by degenerating the ambient Schubert varieties to Hibi toric varieties. We list all possible Calabi--Yau 3-folds of this type up to deformation equivalences, and f\/ind a new example of smooth Calabi--Yau 3-folds of Picard number one; a complete intersection in a locally factorial Schubert variety ${\boldsymbol{\Sigma}}$ of the Cayley plane ${\mathbb{OP}}^2$. We calculate topological invariants and BPS numbers of this Calabi--Yau 3-fold and conjecture that it has a non-trivial Fourier--Mukai partner.}

\Classification{14J32; 14J33; 14M15; 14M25}

\section{Introduction}\label{sec:intro}
Since the discovery of mirror symmetry in the famous example of quintic 3-folds in the projective space, which has a profound relation to conformal f\/ield theories in theoretical phy\-sics~\cite{MR1059831}, constructing examples of Calabi--Yau 3-folds has been playing important roles in the study of mirror symmetry. There are many examples of Calabi--Yau 3-folds of complete intersections in products of (weighted) projective spaces for which we see the symmetry as predicted. In particular, for the examples given as hypersurfaces~\cite{MR1269718} or complete intersections \cite{1993alg.geom.10001B} in Gorenstein toric Fano varieties, the construction of mirror Calabi--Yau 3-folds is now very established based on the combinatorial dualities of ref\/lexive polytopes and nef-partitions. Beyond these combinatorial dualities, so-called the method of conifold transitions gives further examples of mirror sym\-met\-ry for Calabi--Yau 3-folds given as complete intersections in Grassmannians~\cite{MR1619529}, partial f\/lag manifolds \cite{MR1756568}, and also Calabi--Yau 3-folds given as smoothings of hypersurfaces with terminal singularities in Gorenstein toric Fano 4-folds \cite{MR2801412}. The idea of mirror symmetry via conifold transitions is still conjectural, but seems natural since it is based on a kind of monomial-divisor correspondence \cite{MR1253648}, formulated by \cite{MR2112571,MR1673108}.

In this paper, we will consider smooth complete intersection Calabi--Yau 3-folds in the so-called minuscule Schubert varieties (see Def\/inition~\ref{dfn:minuscule}) and study their mirror manifolds via conifold transitions. Classifying such Calabi--Yau 3-folds up to deformation equivalences, we f\/ind a new example of smooth Calabi--Yau 3-fold of Picard number one in a Schubert variety of the Cayley plane $\bO\bP^2$. This example is interesting since there seems to be a non-trivial Fourier--Mukai partner.

Our results for the classif\/ication of smooth complete intersection Calabi--Yau 3-folds in the minuscule Schubert varieties are as follows:
\begin{result}[Propositions \ref{prop:list} and~\ref{prop:invariants}] Let $X$ be a smooth complete intersection Calabi--Yau $3$-fold in a minuscule Schubert variety. Then one of the following holds:
\begin{enumerate}\itemsep=0pt
\item[$(1)$] $X$ is a complete intersection Calabi--Yau $3$-fold of hypersurfaces in homogeneous spaces.
\item[$(2)$] $X$ is a Calabi--Yau $3$-fold of Picard number two.
\item[$(3)$] $X$ is a complete intersection Calabi--Yau $3$-fold of hyperplanes in a Schubert variety ${\boldsymbol{\Sigma}}$ of the Cayley plane $\bO\bP^2$, where ${\boldsymbol{\Sigma}}$ is a~$12$-dimensional locally factorial Fano variety of Picard number one $($Definition~{\rm \ref{dfn:sigma}} and Proposition~{\rm \ref{prop:bSigma})}.
							In this case, the topological invariants of $X$ are
\begin{gather*}
h^{1,1}(X)=1, \qquad h^{2,1}(X)=52,\\ \deg(X)=33,\qquad c_2(X)\cdot H=78,\qquad \chi(X)=-102, \end{gather*}
where $H$ is the ample generator of the Picard group $\Pic(X)\cong \bZ$.
	\end{enumerate}
The cases $(1)$, $(2)$ and $(3)$ include $11$, $2$ and $1$ deformation equivalent classes, respectively.
\end{result}
After establishing this classif\/ication, we will study mirror symmetry for a~Calabi--Yau 3-fold $X$ classif\/ied in Main results~(3). Using the method of conifold transitions, we construct a conjectural mirror family of $X$ and determine the so-called fundamental period following the reference \cite{MR1619529}. We obtain a fourth-order dif\/ferential equation, Picard--Fuchs equation, which annihilate the fundamental period. The fact that the family is mirror symmetric to~$X$ manifests in the monodromy property of this dif\/ferential equation about a distinguished point called a~maximally unipotent monodromy point. Interestingly, we will f\/ind two points which are both maximally unipotent, one of which we identify
with the geometry of~$X$, and the other lead us to conjecture that there exists a~non-trivial Fourier--Mukai partner~$Z$ (Conjecture~\ref{conj:FM}).
As a~supporting evidence of our conjecture, we will obtain integral BPS numbers from the maximally unipotent monodromy points.

\looseness=1 This paper is organized as follows: In Section~\ref{sec:preliminaries} we give some preliminaries. We summarize generalities on posets and distributive lattices, geometry of minuscule Schubert varieties and Hibi toric varieties, and toric degenerations. In particular, we give a combinatorial description of singular locus of minuscule Schubert varieties and Hibi toric varieties, which plays a key role for our proof of Main results. In Section~\ref{sec:list} we make a list of all the deformation equivalent (dif\/feomorphic) classes of smooth complete intersection Calabi--Yau 3-folds in minuscule Schubert varieties. We will f\/ind a new example of such Calabi--Yau 3-folds with Picard number one, embedded in a~locally factorial Schubert variety ${\boldsymbol{\Sigma}}$ in the Cayley plane~$\bO\bP^2$ (Proposition~\ref{prop:list}). In Section~\ref{sec:invariants} we give a computational method of calculating topological invariants for Calabi--Yau 3-folds of Picard number one which degenerate to complete intersection in Gorenstein Hibi toric varieties. We describe the case of ${\boldsymbol{\Sigma}}(1^9)$ in detail as an example, and give a proof of Main results (Proposition~\ref{prop:invariants}). The topological Euler numbers are computed by using a conifold transition.
In Section~\ref{sec:mirror} we study the mirror symmetry of complete intersection Calabi--Yau 3-folds in minuscule Schubert varieties. We construct special one-parameter families of af\/f\/ine complete intersections in~$(\bC^*)^n$ which is conjecturally birational to the mirror families (Conjecture~\ref{conj:mirror}) and give an expression for the fundamental periods (Proposition~\ref{prop:period}). Again, we focus on the example
$X={\boldsymbol{\Sigma}}(1^9)$ and obtain the Picard--Fuchs equation for its mirror, in which we observe two maximally unipotent monodromy points indicating that there is a non-trivial Fourier--Mukai partner of~$X$ (Conjecture~\ref{conj:FM}). We also perform the monodromy calculations and the computations of BPS numbers using mirror symmetry conjecture. Every result seems very similar to that happened to the examples of the Pfaf\/f\/ian--Grassmannian \cite{MR2481271, MR1775415} and the Reye congruence Calabi--Yau 3-fold \cite{MR3212882, MR3166392}. In Appendix~\ref{sec:others}, we discuss other examples in minuscule homogeneous spaces. In Appendix~\ref{sec:bps}, we present the table of BPS numbers calculated for $X={\boldsymbol{\Sigma}}(1^9)$ and its conjectural Fourier--Mukai partner~$Z$.

\subsection*{Note added}
After this work was submitted, there were several developments. Firstly, a dif\/ferent description of the Schubert variety ${\boldsymbol{\Sigma}}$ has been found by \cite{galkin,galkin2}, focusing on an ${\rm SL}(6)$-action on ${\boldsymbol{\Sigma}}$. By using this description, they construct a natural candidate of the Fourier--Mukai partner $Z$ of Calabi--Yau 3-fold $X$ in Conjecture~\ref{conj:FM}. Secondly, it has been found that the Calabi--Yau 3-fold $X$ can be regarded as a complete intersection of a homogeneous vector bundle over a Grassmannian~$G(2,6)$ \cite[Lemma~4.2]{2016arXiv160707821I}. By this description, the $I$-function of $X$ are determined in \cite{2016arXiv160708137I}. Physics related to~$X$ has also been investigated in \cite{2015arXiv150500099G,2016arXiv161206214H,2016arXiv160202487U}.

\section{Preliminaries}\label{sec:preliminaries}

\subsection{Posets and distributive lattices}\label{subsec:posets}

A partially ordered set is called a \textit{poset} for short. Let $u$ and $v$ be distinct elements of a~poset~$P$, then there are three possibilities, i.e., $u \prec v$, $v \prec u$ or $u$ and $v$ are incomparable, written as $u \not\sim v$. We say that $u$ \textit{covers} $v$ if $u \succ v$ and there is no $w\in P$ with $u \succ w \succ v$. The \textit{Hasse diagram} of a poset $P$ is the oriented graph with the vertex set~$P$, having an edge $e=\{u, v\}$ going down from $u$ to $v$ whenever $u$ covers $v$ in $P$. We can read a lot of information for a poset by drawing the Hasse diagram in a plane. For a poset $P$, an \textit{order ideal} is a subset $I \subset P$ with the property that
\begin{gather*}	u\in I \quad \text{and} \quad v\prec u \quad \text{imply} \quad v\in I.\end{gather*}

A \textit{lattice} $L$ is a poset for which each pair of elements $\alpha, \beta \in L$ has a least upper bound $\alpha \vee \beta$ (called the join) and a greatest lower bound $\alpha \wedge \beta$ (called the meet) in $L$. A \textit{distributive lattice} is a lattice on which the following identity holds for all triple elements $\alpha, \beta, \gamma \in L$
\begin{gather*}
\alpha\wedge (\beta\vee \gamma) = (\alpha\wedge \beta)\vee (\alpha\wedge \gamma).
\end{gather*}
It is easy to see that a f\/inite lattice $L$ has a unique maximal and minimal element with respect to the partial order on~$L$. In a lattice~$L$, an element $\alpha$ is said to be \textit{join irreducible} if $\alpha$ is neither the minimal element nor the join of a f\/inite set of other elements.

For a f\/inite poset $P$, the order ideals of $P$ form a distributive lattice $J(P)$ with the partial order given by set inclusions. The join and the meet on $J(P)$ correspond to the set union and the set intersection, respectively. On the other hand, the full subposet of join irreducible elements of $J(P)$ coincides with $P$ by the Birkhof\/f representation theorem \cite{MR598630}. In fact, this gives a one-to-one correspondence between f\/inite posets and f\/inite distributive lattices. We give an example of this correspondence in the following.

\begin{Example}\label{ex:poset}
We introduce an example which appears repeatedly throughout this paper. Let us def\/ine a poset $\bfP$ by its Hasse diagram drawn on the left side of Fig.~\ref{fig:hasse}.
\begin{figure}[h!]\centering
\includegraphics[scale=0.8]{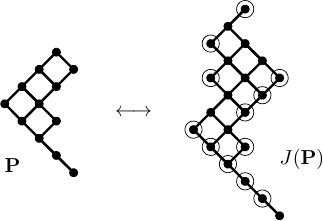}
\caption{The Birkhof\/f correspondence.}\label{fig:hasse}
\end{figure}
The distributive lattice $J(\bfP)$ is obtained by drawing its Hasse diagram (on the right side of Fig.~\ref{fig:hasse}) while focusing on the inclusion relation of order ideals of~$\bfP$ explicitly, which is regarded as the partial order on~$J(\bfP)$. For instance, the upper part of the Hasse diagram of~$J(\bfP)$ in Fig.~\ref{fig:hasse} corresponds to the following order ideals of~$\bfP$:
\begin{figure}[h!]\centering
\includegraphics[scale=0.8]{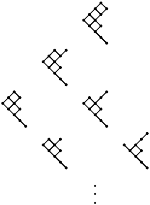}
\end{figure}

Note that the maximal element and the minimal element of $J(\bfP)$ correspond to the whole $\bfP$ and the empty set $\varnothing$, respectively. The circled vertices of $J(\bfP)$ in Fig.~\ref{fig:hasse} represent join irreducible elements of $J(\bfP)$, i.e., the vertices with exactly one edge below. We can reconstruct the poset~$\bfP$ as the set of circled vertices with the induced order in~$J(\bfP)$.
\end{Example}

Now we introduce further def\/initions on a f\/inite poset. Def\/ine the source $s(e) =u$ and the target $t(e)=v$ for an edge $e=\{u, v\}$ of the Hasse diagram of $P$ if $u$ covers~$v$. A \textit{chain} of length~$k$ in $P$ is a sequence of elements $u_0 \prec u_1 \prec \cdots \prec u_k \in P$. A chain is called \textit{maximal} if there is no $v \prec u_0$ or $w \succ u_k$ in $P$ and $u_i$ covers $u_{i-1}$ for all $1\le i \le k$. We call a f\/inite poset $P$ \textit{pure} if every maximal chain has the same length. Let us def\/ine the associated bounded poset $\hat{P}=P\cup \{ \hat{0}, \hat{1} \}$ for any f\/inite poset~$P$, by extending the partial order on $P$ with $\hat{0} \prec u \prec \hat{1}$ for all $u\in P$. Every f\/inite poset $P$ has a height function $h$ by def\/ining $h(u)$ to be the length of the longest chain bounded above by $u$ in $\hat{P}$. We also def\/ine the \textit{height} $h_P$ of $P$ as $h(\hat{1})$. For example, $h_\bfP=9$ for the pure poset~$\bfP$ in Example~\ref{ex:poset}.

\subsection{Minuscule Schubert varieties}\label{subsec:Schubert}
We give a brief summary on minuscule homogeneous spaces and minuscule Schubert varieties. A basic reference of this subject is~\cite{LMS1979}. Let~$G$ be a~simply connected simple complex algebraic group, $B$ a Borel subgroup and $T$ a maximal torus in~$B$. We denote by $R^+$ the set of positive roots and
by $S=\{ \alpha_1,\dots,\alpha_n \}$ the set of simple roots. Let $W$ be the Weyl group of~$G$. Denote by~$\Lambda$ the character group of~$T$, which is also called the weight lattice of~$G$. The weight lattice~$\Lambda$ is generated by the fundamental weights $\lambda_1, \dots, \lambda_n$ def\/ined by $(\alpha_i^\vee , \lambda_j )=\delta_{ij}$ for $1\le i, j \le n$, where~$(\,,\,)$ is a $W$-invariant inner product and $\alpha^\vee := 2 \alpha/(\alpha,\alpha)$. An integral weight \mbox{$\lambda=\sum n_i \lambda_i \in\Lambda$} is said to be dominant if $n_i \ge 0$ for all $i=1,\dots , n$. For an integral dominant weight $\lambda\in \Lambda$, we denote by $V_\lambda$ the irreducible $G$-module of highest weight~$\lambda$. The associated homogeneous space~$G/Q$ of~$\lambda$ is the $G$-orbit of the highest weight vector in the projective space~$\bP(V_\lambda)$, where $Q\supset B$ is the associated parabolic subgroup of~$G$. A~Schubert variety in~$G/Q$ is the closure of a~$B$-orbit in~$G/Q$. Recall the following def\/inition of a~minuscule weight, a minuscule homogeneous space and minuscule Schubert varieties.

\begin{Definition}[cf.\ {\cite[Def\/inition 2.1]{LMS1979}}] \label{dfn:minuscule} Let $\lambda \in \Lambda$ be a fundamental weight. We call $\lambda$ \textit{minuscule} if it satisf\/ies the following equivalent conditions:
\begin{enumerate}\itemsep=0pt
\item[(1)] 
Every weight of $V_\lambda$ is in the orbit $W\lambda \subset \Lambda$.
\item[(2)] 
$( \alpha^\vee , \lambda ) \leq 1$ for all $\alpha \in R^+$.
\end{enumerate}
The homogeneous space $G/Q$ associated with a minuscule weight $\lambda$ is said to be minuscule. The Schubert varieties in minuscule $G/Q$ are also called minuscule.
\end{Definition}

Let us recall that a parabolic subgroup $Q \supset B$ is determined by a subset $S_Q$ of $S$ asso\-cia\-ted with negative root subgroups. A useful notation for a homogeneous space $G/Q$ is to cross the nodes in the Dynkin diagram which correspond to the simple roots in $S\setminus S_Q$. With this notation, the minuscule homogeneous spaces are as shown in Table~\ref{tbl:minuscule}. This contains the Grassmannians~$G(k,n)$, the orthogonal Grassmannians~$OG(n,2n)$, even dimensional quadrics~$Q^{2n}$ and, f\/inally, the Cayley plane $\bO \bP^2=E_6/Q_1$ and the Freudenthal variety $E_7/Q_7$, where we use the Bourbaki labelling for the roots. We omit minuscule weights for groups of type $B$ and type~$C$, since they give isomorphic varieties to those for simply laced groups.
\begin{table}[h]\centering
 \includegraphics{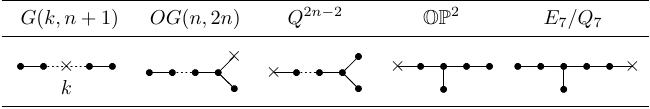}
 \caption{Minuscule homogeneous spaces.}\label{tbl:minuscule}
\end{table}

The Weyl group $W$ is generated by simple ref\/lections $s_\alpha \in W$ for $\alpha \in S$. These generators def\/ine the length function $l$ on $W$. Let us denote by $W_Q$ the Weyl group of $Q$, i.e., the subgroup generated by $\{s_\alpha \in W \,|\, \alpha \in S_Q \}$, and by $W^Q$ the set of minimal length representatives of the coset $W/W_Q$ in $W$. For any $w \in W^Q$, we denote by $X(w)=\overline{BwQ/Q}$ the Schubert variety in $G/Q$ associated with $w$, which is an $l(w)$-dimensional normal Cohen--Macaulay projective variety with at worst rational singularities. There is a natural partial order $\prec$ on $W^Q$
called the Bruhat order, def\/ined by $w_1 \preceq w_2 \Leftrightarrow X(w_1)\subset X(w_2)$. We recall the following fundamental fact for minuscule homogeneous spaces.

\begin{Proposition}[{\cite[Proposition V.2]{MR2941004}}] \label{prop:distr}
For a minuscule homogeneous space $G/Q$, the po\-set~$W^Q$ is a finite distributive lattice.
\end{Proposition}

From Proposition \ref{prop:distr} and the Birkhof\/f representation theorem, we can def\/ine the \textit{minuscule poset} $P_Q$ for the minuscule $G/Q$ such that $J(P_Q)=W^Q$ as in~\cite{MR2941004}. Moreover, the order ideal $P_w \subset P_Q$ associated with $w\in W^Q$ is called the minuscule poset for the minuscule Schubert variety $X(w)\subset G/Q$. For example, the order ideals $P_Q$, $\varnothing \subset P_Q$ turn out to be the minuscule posets for the total space $X(w_Q)=G/Q$ and the $B$-f\/ixed point $X(\id)=Q/Q$, respectively, where $w_Q$ is the unique longest element in $W^Q$. The minuscule poset for minuscule Schubert varieties is a generalization of the Young diagram for Grassmann Schubert varieties.
\begin{Example} \label{ex:sigma} An easy method to compute the Hasse diagram of $W^Q$ is to trace out the $W$-orbit of certain dominant weight whose stabilizer coincides with $W_Q$ (cf.\ \cite[Section~4.3]{MR1038279}). Denote by $(i j \cdots k)$ the element $w=s_{\alpha_i} s_{\alpha_j} \cdots s_{\alpha_k} \in W$, where $s_\alpha$ is the simple ref\/lection with respect to $\alpha \in S$. The initial part of the Hasse diagram of $W^Q$ for the Cayley plane $\bO\bP^2$ is the following
\begin{table}[h]\centering
 \includegraphics{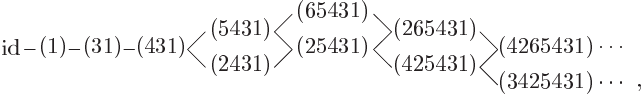}
\end{table}

\noindent
where the right covers the left for connected two elements with respect to the Bruhat order. Thus we obtain the Hasse diagram of the distributive lattice~$W^Q$ and hence the minuscule poset for every Schubert variety in~$\bO\bP^2$.
\end{Example}
\begin{Definition} \label{dfn:sigma} In the above notation, let us set $w=(314354265431)\in W^Q$. We denote by ${\boldsymbol{\Sigma}}$ the associated $12$-dimensional Schubert variety $X(w)$ in the Cayley plane $\bO\bP^2$, which corresponds to the minuscule poset $\bfP$ in Example~\ref{ex:poset}.
\end{Definition}

\begin{Remark} \label{rem:minuscule} We remark that in addition to	(1) and (2) of Def\/inition~\ref{dfn:minuscule} there is a~geometric characterization of minuscule homogeneous spaces in \cite[Def\/inition~2.1]{LMS1979}. A~fundamental weight~$\lambda$ is minuscule if and only if the following condition holds.
 \begin{enumerate}\itemsep=0pt
\item[(3)]	For the associated homogeneous space $G/Q$, the Chevalley formula
	\begin{gather*} [H]\cdot [X(w)] =\sum_{w \text{ covers } w'} [X(w')] \end{gather*}
	holds for all $w\in W^Q$ in the Chow ring of $G/Q$, where $H$ is the unique Schubert divisor in $G/Q$ and	$W^Q$ is the poset with the Bruhat order.
 \end{enumerate}

 As a corollary, it turns out that the degree of
 a minuscule Schubert variety $X(w)$
 with respect to $\scO_{G/Q}(1)|_{X(w)}$
 equals the number of maximal chains in $J(P_w)$.
 For example, we obtain $\deg {\boldsymbol{\Sigma}}=33$ by counting the maximal chains
 in $J(\bfP)$ in Fig.~\ref{fig:hasse}.
\end{Remark}

We introduce further def\/initions in order to describe singularities of minuscule Schubert varieties. As we expect from the computation in Example~\ref{ex:sigma}, the Bruhat order on $W^Q$ is generated by simple ref\/lections for minuscule $G/Q$ \cite[Lemma~1.14]{MR1080976}, that is,
\begin{gather}				\label{eq:fact}
				w_1 \text{ covers } w_2 \Leftrightarrow w_1=s_\alpha \cdot w_2
				\text{ and } l(w_1)=l(w_2)+1 \text{ for some }\alpha \in S. \end{gather}
From this fact, a join irreducible element $u\in W^Q$ covers the unique element \smash{$s_{\beta_Q(u)}\cdot u\in W^Q$} where $\beta_Q(u)\in S$. Thus we can def\/ine a natural map $\beta_Q\colon P_Q \rightarrow S$. We call $\beta_Q$ the \textit{natural co\-lo\-ra\-tion} for a minuscule poset $P_Q$ by simple roots $S$. We also def\/ine the natural coloration~$\beta_w$ on each minuscule poset $P_w \subset P_Q$ by restricting~$\beta_Q$ on~$P_w$. The minuscule poset~$P_w$ with the natural coloration $\beta_w\colon P_w \rightarrow S$ has in fact the same information as the minuscule quiver introduced by Perrin \cite{MR2360316, MR2466424}, which gives a good description of geometric properties of minuscule Schubert varieties $X(w)$. We translate some combinatorial notion and useful facts on the geometry of minuscule Schubert varieties $X(w)$ from the reference \cite{MR2360316, MR2466424} into our terminology.

\begin{Definition} \label{dfn:peaksandholes} Let $P$ be a minuscule poset with the natural coloration $\beta\colon P\rightarrow S$.
 \begin{enumerate}\itemsep=0pt
		 \item[(1)] A \textit{peak} of $P$ is a maximal element $u$ in $P$.
		 \item[(2)] A \textit{hole} of $P$ is a maximal element $u$
			in $\beta^{-1}(\alpha)$ for some $\alpha \in S$
			such that there are exactly two elements $v_1, v_2\in P$
 			with $u \prec v_i$ and
			$( \beta(u)^\vee , \beta(v_i) ) \ne 0$ ($i=1, 2$).
 \end{enumerate}
 We denote by $\mathrm{Peaks}(P)$ and $\mathrm{Holes}(P)$
 the set of peaks and holes of $P$, respectively.
 A hole $u$ of the poset $P$ is said to be \textit{essential}
 if the order ideal
 $P^{u}:= \left\{ v\in P \,|\, v\not\succcurlyeq u \right\}$ contains
 all other holes in $P$.
\end{Definition}

Let $X(w)$ be a minuscule Schubert variety in $G/Q$ and $P_w$ the associated minuscule poset. Weil and Cartier divisors on $X(w)$ are described in terms of the poset $P_w$. In fact, it is clear that any Schubert divisor coincides with a Schubert variety $D_u$ associated with $P_w^u$ for some $u \in \mathrm{Peaks}(P_w)$. It is well-known that the divisor class group $\mathrm{Cl}(X(w))$ is the free $\bZ$-module generated by the classes of the Schubert divisors $D_u$ for $u\in \mathrm{Peaks}(P_w)$, and the Picard group $\mathrm{Pic}(X(w))$ is isomorphic to $\bZ$ generated by $\scO_{G/Q}(1)|_{X(w)}$. As we saw in Remark~\ref{rem:minuscule}, the Cartier divisor corresponding to $\scO_{G/Q}(1)|_{X(w)}$ is given by
\begin{gather*}\sum_{u\in \mathrm{Peaks}(P_w)} D_u.\end{gather*}

We use the following results by Perrin about the singularities of $X(w)$.

\begin{Proposition}[{\cite{MR2360316, MR2466424}}] \label{prop:schubertsing}
 Let $X(w)$ be a minuscule Schubert variety and $P_w$ the associated minuscule poset.
 \begin{enumerate}\itemsep=0pt
	\item[$(1)$] An anticanonical Weil divisor of $X(w)$ is given by
		\begin{gather*} -K_{X(w)}= \sum_{u\in \mathrm{Peaks}(P_w)} (h(u)+1) D_u.\end{gather*}
	 In particular, $X(w)$ is Gorenstein if and only if~$P_w$ is pure. In this case $X(w)$ is a Fano variety of index $h_{P_w}$
		{\rm (}see {\rm \cite[Proposition~4.17]{MR2360316})}.
	\item[$(2)$]
	 The Schubert subvariety associated with the order ideal $P_w^{u}\subset P_w$ for an essential hole~$u$ of~$P_w$ is an irreducible component
	 of the singular locus of~$X(w)$. All the irreducible components of the singular locus are obtained in this way $($see {\rm \cite[Theorem~2.7(1)]{MR2466424})}.
 \end{enumerate}
\end{Proposition}

We apply Proposition \ref{prop:schubertsing} to our example ${\boldsymbol{\Sigma}}$ and obtain the following proposition.

\begin{Proposition}\label{prop:bSigma}
 Let ${\boldsymbol{\Sigma}}$ be the minuscule Schubert variety in $\bO\bP^2$ defined in Definition~{\rm \ref{dfn:sigma}}.
 \begin{enumerate}\itemsep=0pt
	\item[$(1)$] ${\boldsymbol{\Sigma}}$ is a locally factorial Gorenstein Fano variety of index $9$.
	\item[$(2)$] The singular locus of ${\boldsymbol{\Sigma}}$ is isomorphic to $\bP^5$.
 \end{enumerate}
\end{Proposition}
\begin{proof} The former holds because the corresponding minuscule poset $\bfP$ (Fig.~\ref{fig:hasse}) is pure with $h_{\bfP}=9$ and has a unique peak. From~(\ref{eq:fact}) and the computation of the Hasse diagram of $W^Q$ for $\bO\bP^2$ in Example~\ref{ex:sigma}, the natural coloration $\beta \colon \bfP \rightarrow S$ can be determined. The resulting coloration $\beta$ is represented as a vertical projection in the right in Fig.~\ref{fig:color}, where we identify the simple roots $S$ and vertices of the Dynkin diagram of $E_6$.
\begin{figure}[h]\centering
 \includegraphics{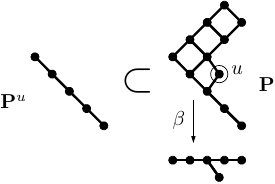}

 \caption{The natural coloration of the minuscule poset $P$ and the singular locus.}
 \label{fig:color}
\end{figure}
In Fig.~\ref{fig:color}, the unique (essential) hole of $\bfP$ is the circled vertex $u$, and its color is read as $\alpha_2\in S$. The corresponding Schubert subvariety is described by the minuscule poset $\bfP^{u}$, which coincides with the singular locus of ${\boldsymbol{\Sigma}}$ by Proposition~\ref{prop:schubertsing}. It is isomorphic to $\bP^5$ because the degree equals one.
\end{proof}

Finally, we recall the vanishing theorems for minuscule Schubert varieties for later use.
\begin{Theorem}[{\cite[Theorem~7.1]{LMS1979}}] \label{thm:vanishing} Let $\lambda$ be a minuscule weight, $G/Q \subset \bP(V_\lambda)$ the asso\-cia\-ted homogeneous space and $X(w)\subset G/Q$ a minuscule Schubert variety.
	Then the following properties hold:
	 \begin{enumerate}\itemsep=0pt
		\item[$(1)$] $H^0(\bP(V_\lambda),\scO(m))\rightarrow H^0(X(w),\scO(m))$ is surjective for all $m \ge 0$,
		\item[$(2)$] $H^i(X(w),\scO(m))=0$ for all $m \in \bZ$ and $0 < i< l(w)$,
		\item[$(3)$] $H^{l(w)}(X(w),\scO(m))=0$ for all $m \ge 0$.
	 \end{enumerate}
\end{Theorem}

\subsection{Hibi toric varieties}

We introduce projective varieties called Hibi toric varieties. The geometry of minuscule Schubert varieties will be related with the corresponding geometry of Hibi toric varieties.

Let us start with def\/ining the polytope $\Delta(P)$ called the \textit{order polytope} for a f\/inite poset $P$. Order polytopes has been studied in many literatures of combinatorics, for example~\cite{MR824105}. Let $N=\bZ P$ and $M=\Hom_\bZ(N,\bZ)$ be free $\bZ$-modules of rank $\abs P$, and $N_\bR$, $M_\bR$ the real scalar extensions $N\otimes_\bZ \bR$, $M\otimes_\bZ \bR$, respectively. We def\/ine the order polytope $\Delta(P)\subset M_\bR$ by
\begin{gather*}
\Delta(P):=\big\{ (x_u)_{u\in P} \, |\, 0 \leq x_u \leq x_v \leq 1 \text{ for all } u \prec v \in P \big\}.
\end{gather*}
It is easy to see that the order polytope $\Delta(P)$ is an integral convex polytope of dimension $\abs P$.
\begin{Definition}
Let $P$ be a f\/inite poset and $\Delta(P)$ the order polytope for $P$. The projective toric variety associated with $\Delta(P)$,
\begin{gather*}
				\bP_{\Delta(P)}:= \Proj \bC[\Cone\left(1\times \Delta(P)\right) \cap (\bZ \times M)]
\end{gather*}
is called the \textit{$($projective$)$ Hibi toric variety} for $P$.
\end{Definition}

\begin{Example}One of the simplest examples of order polytopes is the Gelfand--Tsetlin polytopes for fundamental weights of special linear groups ${\rm SL}(n+1,\bC)$. In general, the Gelfand--Tsetlin polytope for an integral dominant weight $\lambda=(\Lambda_0,\dots,\Lambda_n)\in
\bZ^{n+1}/\langle(1,\dots,1)\rangle \simeq \Lambda$ is def\/ined by the following inequalities in $\bR^{n(n+1)/2}$
\begin{gather*}
\Lambda_n \le x_{i+1,j+1}\le x_{i,j}\le x_{i,j+1} \le \Lambda_i \quad \text{for all }0\le i\le j\le n-1.
\end{gather*}
It is standard to represent these in a diagram like the left of Fig.~\ref{fig:gt} for $n=3$.
\begin{figure}[h]\centering
 \includegraphics{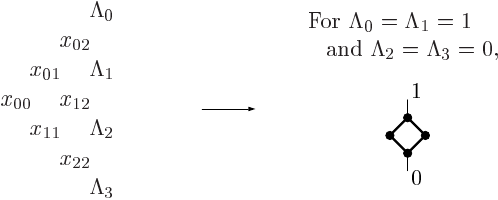}
 \caption{A Gelfand--Tsetlin polytope for ${\rm SL}(4,\bC)$.} \label{fig:gt}
\end{figure}
The Gelfand--Tsetlin polytope for a fundamental weight $\lambda=(1,\dots , 1, 0, \dots ,0)$ is the order polytope for a~poset of rectangle shape, like the right of Fig.~\ref{fig:gt} for example. It coincides with the minuscule poset for a Grassmannian of type $A$. The corresponding Hibi toric variety $\bP_{\Delta(P)}$ is the toric variety $P(k,n+1)$ def\/ined by~\cite{MR1619529}.
\end{Example}

\begin{Remark} \label{rem:alternative}
There is an alternative def\/inition of a Hibi toric variety using its homogeneous coordinate ring. Let $J(P)$ be the distributive lattice of order ideals of a f\/inite poset~$P$. Denote by $\bC[J(P)]$ the polynomial ring over $\bC$ in $\abs{J(P)}$ indeterminates $p_\alpha$ $(\alpha\in J(P))$. Let $I(J(P))\subset \bC[J(P)]$ the homogeneous ideal generated by the following binomial relations:
\begin{gather*}
	p_\tau p_\phi - p_{\tau\wedge\phi} p_{\tau\vee\phi} \quad ( \tau \not\sim \phi).
\end{gather*}
One can show that the graded algebra $A_{J(P)}:=\bC [J(P) ]/I(J(P))$ with the standard $\bZ$-grading inherited from $\bC [J(P) ]$ coincides with the homogeneous coordinate ring of the Hibi toric va\-riety~$\bP_{\Delta(P)}$ with the embedding associated with $\Delta(P)$. The graded algebra $A_{J(P)}$ is usually called \textit{Hibi algebra} on the distributive lattice $J(P)$ (cf.~\cite{MR951198}).
\end{Remark}

We will use the following facts known for the Hibi algebra introduced above.
\begin{Proposition}[{\cite[Corollary of Lemma 5]{MR790025}}]\label{prop:deg}Let $P$ be a finite poset and $J(P)$ the distributive lattice of order ideals of $P$. The Poincar\' e series of the Hibi algebra $A_{J(P)}$ has the following form:
	\begin{gather*} \scP(t)= \sum_{i=0}^{\abs P} c_i (t/(1-t))^i, \end{gather*}
where $c_0 = 1$ and $c_i$ $(i\ne 0)$ is the number of chains of length $i$ in $J(P)$. In particular, the degree of $\bP_{\Delta(P)} \subset \Proj \bC[J(P)]$ equals the number of maximal chains in $J(P)$.
\end{Proposition}

\begin{Proposition}[{\cite{MR951198}}]\label{prop:hibi} Let $P$ be a finite poset. The Hibi algebra $A_{J(P)}$ is Gorenstein if and only if $P$ is pure.
\end{Proposition}

One of the nice properties of Hibi toric varieties is that torus invariant subvarieties are also Hibi toric varieties. We need explicit descriptions of these subvarieties and their singularities. We summarize the results obtained in \cite{MR1382045} and also in \cite{MR2416742, MR2677619} and \cite{MR2761126}.

Let us recall the inequalities in the def\/inition of an order polytope $\Delta(P)$. They are generated by $x_{s(e)} \geq x_{t(e)}$ for all $e \in E$, where $E$ is the set of edges of the Hasse diagram of $\hat{P}=P\cup \left\{ \hat0, \hat1 \right\}$ and $x_{\hat 0}=0$ and $x_{\hat 1}=1$. We obtain a face of $\Delta(P)$ by replacing some of these inequalities with equalities. To describe the faces of~$\Delta(P)$, we move def\/initions about posets.

Recall that a \textit{full subposet} $y \subset \hat P$ to be a~subset of $\hat P$ whose poset structure is that inherited from $\hat P$. We call a full subposet $y\subset \hat P$ \textit{connected} if all the elements in $y$ are connected by edges in the Hasse diagram of $y$, \textit{convex} if $u, v \in y$ and $u \prec w \prec v$ imply $w \in y$, and a \textit{cycle} if the Hasse diagram of $y$ is a cycle as an unoriented graph (e.g., the poset of the right of Fig.~\ref{fig:gt} is a~cycle).
\begin{Definition}\label{dfn:admissible}A surjective map $f\colon \hat P \rightarrow \hat P'$ is called a \textit{contraction} of the bounded poset $\hat P$ if every f\/iber $f^{-1}(i)$ $(i \in \hat P ')$ is a connected full subposet of $\hat P$ not containing both $\hat 0$ and $\hat 1$, and 	the following condition holds for all $u_k, v_k \in f^{-1}(k)$ $(k \in \hat P')$
	and $i \neq j\in \hat P'$:
	\begin{gather*}
					\text{a relation }u_i \prec u_j \text{ implies } v_i \not\succ v_j.\end{gather*}
 \end{Definition}
\begin{Remark}
A contraction $f\colon \hat P \rightarrow \hat P'$ gives a natural partial order on the image set $\hat{P'}$, i.e., the partial order generated by the following relations:
	\begin{gather*}
					i \prec j \Leftrightarrow \text{there exist }u \in f^{-1}(i) \text{ and }v \in f^{-1}(j)
					\text{ such that }u \prec v \text{ in }\hat{P}.\end{gather*}
Further, $\hat{P'}$ becomes a bounded poset by setting $\hat{1}\in f^{-1}(\hat{1})$ and $\hat{0}\in f^{-1}(\hat{0})$. In fact, the above def\/inition of contraction coincides with the more abstract def\/inition called the f\/iber-connected tight surjective morphism of bounded posets (see~\cite{MR1382045}).	
\end{Remark}
For a contraction $f\colon \hat{P} \rightarrow \hat{P'}$,
the corresponding face of $\Delta(P)$ is given by
\begin{gather*}
\theta_{f} := \big\{ x \in \Delta(P)\, |\, x_u=x_v \text{ for all }u, v \in f^{-1}(i) \text{ and }i \in \hat{P'}\big\}.
\end{gather*}
Conversely, we can reconstruct the contraction from each face
$\theta_f \subset \Delta(P)$
by looking at the coordinates of general point in $\theta_f$.
Now we can summarize the classical fact on the face structure of
order polytopes in our terminology (cf.\ \cite[Theorem~1.2]{MR1382045}).

\begin{Proposition}\label{prop:facestr}
 Let $P$ be a finite poset, and $\Delta(P)$ the associated order polytope.
 The above construction gives
 a one-to-one correspondence between faces of $\Delta(P)$ and
 contractions of $\hat{P}$.
 Moreover, an inclusion of the faces
 corresponds to a composition of contractions.
\end{Proposition}

\begin{Remark}\label{rem:hibi}\quad\samepage
\begin{enumerate}\itemsep=0pt
\item[(1)] As a special case, there is a one-to-one correspondence between facets and edges $E$ of the Hasse diagram of $\hat P$. For each $e\in E$, the corresponding contraction $f_e\colon \hat P\rightarrow \hat P'$, and hence the facet $\theta_{f_e}$, is easily obtained by replacing the inequality $x_{s(e)} \ge x_{t(e)}$ by the equality $x_{s(e)} = x_{t(e)}$, where $x_{\hat 0}=0$ and $x_{\hat 1}=1$. We also write $\theta_e=\theta_{f_e}$ to represent the facet corresponding to $e\in E$.
\item[(2)] It is obvious that the face $\theta_{\hat{P} \rightarrow \hat{P'}}\subset \Delta(P)$ coincides with the $\abs{P'}$-dimensional order polytope $\Delta(P')$ under a suitable choice of subspace of $M_\bR$ and a unimodular transformation. This means that the torus invariant subvarieties in Hibi toric varieties are also Hibi toric varieties as noted above.
\end{enumerate}
\end{Remark}

We note the following nice description of singularities of Hibi toric varieties.
\begin{Theorem}[{\cite[Theorem 2.3 and proof of Corollary 2.4]{MR1382045}}]	\label{thm:sing}
Let $\bP_{\Delta(P)}$ be a Hibi toric variety for a finite poset $P$. A face $\theta_f\subset \Delta(P)$ corresponds to an irreducible component of the singular loci of $\bP_{\Delta(P)}$ if and only if one of the fiber $f^{-1}(i)$ of the contraction $f\colon \hat P \rightarrow \hat P'$ is a~minimal convex cycle and all other fibers $f^{-1}(j)$ $(j \ne i)$ consist of one element.
\end{Theorem}

As in the case of minuscule Schubert varieties, Weil and Cartier divisors on a Hibi toric variety $\bP_{\Delta(P)}$ are naturally described in terms of the poset $P$. One can show that the divisor class group $\mathrm{Cl}(\bP_{\Delta(P)})$ is a free $\bZ$-module of rank $\abs E -\abs P$ and the Picard group $\Pic(\bP_{\Delta(P)})$ is a free $\bZ$-module of rank $c$, where $c$ is the number of connected components of the Hasse diagram of $P$. More precisely, $\mathrm{Cl}(\bP_{\Delta(P)})$ is presented by toric divisors modulo the following linear equivalence relations,
\begin{gather*}
\sum_{s(e)=u} D_e \simeq \sum_{t(e)=u} D_e, \qquad u\in P,
\end{gather*}
where $D_e$ denotes the toric divisor on $\bP_{\Delta(P)}$ corresponding to $e \in E$. Similarly, $\Pic(\bP_{\Delta(P)})$ is generated by the invertible sheaves associated with
\begin{gather*}
		\sum_{\substack{s(e)=\hat 1\\t(e)\in y}}D_e \simeq 	\sum_{\substack{s(e)\in y \\t(e)= \hat 0}}D_e
\end{gather*}
for all connected components $y \subset P$. These properties follow from almost the same arguments done for a special case \cite[Lemma~3.2.2]{MR1756568} based on the standard technologies on toric geometry.

If the f\/inite poset $P$ is pure, the Hibi toric variety $\bP_{\Delta(P)}$ becomes a Gorenstein Fano variety of index $h_P$. In fact, the anticanonical divisor turns out to be an ample Cartier divisor,
\begin{gather*}
				-K_{\bP_{\Delta(P)}}=\sum_{e\in E}D_e= \sum_{k=1}^{h_P} D_{E^k}	\simeq h_P D_{E^1},
\end{gather*}
where $D_{E'}:= \sum\limits_{e\in E'} D_e$ for any subset $E'\subset E$ and $E^k:= \{ e\in E \,|\, h(s(e))=k \}$. Note that $\scO(D_{E^k})$ coincides with a very ample invertible sheaf $\scO(1)$ associated with the polytope $\Delta(P)$. In this case, it is known that the singularities of $\bP_{\Delta(P)}$ are rather mild as in the following lemma:

\begin{Lemma}[{\cite[Lemma 1.4]{MR2770546}}]\label{lem:hibi2} Let $P$ be a finite pure poset. The Hibi toric variety $\bP_{\Delta(P)}$ has at worst terminal singularities.
\end{Lemma}

We summarize some more properties on Gorenstein Hibi toric varieties and complete intersection Calabi--Yau manifolds of hypersurfaces in it.

Let $P$ be a f\/inite pure poset. The Hibi toric variety $\bP_{\Delta(P)}$ is a Gorenstein Fano variety of index $h_P$ as we saw. By \cite[Theorem~4.1.9]{MR1269718}, a unimodular transformation of $h_P\Delta(P)$ becomes a \textit{reflexive polytope} corresponding to the anticanonical sheaf $-\scK_{\bP_{\Delta(P)}}=\scO(h_P)$, i.e., it contains the unique internal integral point $0$ and every facet has integral distance one from $0$. We take $\Delta:= \sum\limits_{u\in P}h(u)\chi_u- h_P \Delta(P)$ as such a ref\/lexive polytope, where $\left\{ \chi_u \,|\, u\in P \right\}\subset M$ is the dual basis of $P\subset N$, i.e., $\chi_u(v)=\delta_{uv}$ for $u,v\in P$. Denote by $l(\theta)=\abs{\theta \cap M}$ the number of integral points in a polytope $\theta \subset M_\bR$, and $l^*(\theta)$ the number of integral points in the interior of $\theta$. We note that $l(\Delta(P))=\abs{J(P)}$ and $l^*(h_P \Delta(P))=l^*(\Delta)=1$.

A \textit{maximal projective crepant partial resolution $($MPCP-resolution$)$} of $\bP_{\Delta(P)}$ is a crepant morphism $\hat \bP_{\Delta(P)} \rightarrow \bP_{\Delta(P)}$ where $\hat \bP_{\Delta(P)}$ has only $\bQ$-factorial terminal singularities. By \cite[Theo\-rem~2.2.24]{MR1269718}, there exists at least one MPCP-resolution for a Gorenstein toric Fano variety. The MPCP-resolution of $\bP_{\Delta(P)}$ is obtained by maximal triangulations on~$N$ of each facet of
the polar dual polytope $\Delta^*:= \big\{t\in N_\bR \,|\, x(t) \ge -1 \text{ for all } x \in \Delta \big\}$ of $\Delta$.

Let $X_0 \subset \bP_{\Delta(P)}$ be a complete intersection of $r$ general hypersurfaces of degree $(d_1, \dots , d_r)$ where $\sum\limits_{j=1}^r d_j=h_P$.
Denote by $Y=\hat{X_0}$ the strict transform of $X_0$ in a MPCP-resolution $\hat{\bP}_{\Delta(P)}\rightarrow \bP_{\Delta(P)}$. We have convenient formulas for some of the stringy Hodge numbers of $X_0$ in Theorem~\ref{thm:hodge}. The stringy Hodge numbers of a singular variety $X_0$ coincide with the usual Hodge numbers of $Y$ if there exists a~crepant resolution $Y\rightarrow X_0$.

\begin{Theorem} \label{thm:hodge}
The stringy $(1,*)$-Hodge numbers of a general complete intersections $X_0$ of degree $(d_1, \dots, d_r)$ with $\sum\limits_{j=1}^r d_j=h_P$ in a Gorenstein Hibi toric variety $\bP_{\Delta(P)}$ are given by
 \begin{gather}	\label{eq:formula}
 h_{\mathrm{st}}^{1,1}(X_0)=\abs E- \abs P,\qquad h_{\mathrm{st}}^{1,k}(X_0)=0, \qquad 1<k<\abs P- r-1,\\
 h_{\mathrm{st}}^{1,\abs P -r-1}(X_0)= \sum_{i\in I}\left[ \sum_{J \subset I}(-1)^{\abs J}	l\left((d_i-d_J)\Delta(P) \right)\right]-
 \sum_{J \subset I}(-1)^{r-\abs J}\left[ \sum_{e\in E}	l^*(d_J\theta_e) \right] -\abs P,\nonumber
\end{gather}
where $I=\{1, \dots, r\}$, $d_J=\sum\limits_{j\in J}d_j$ and $\theta_e$ is the facet of $\Delta(P)$ corresponding to an edge \smash{$e\in E$}.	A~nonzero contribution in the first term of $h_{\mathrm{st}}^{1,\abs P -r-1}(X_0)$ comes only from the range of \smash{$d_i-d_J\ge 0$} and in the second term from that of $d_J = h_P-1$ or $h_P$.
\end{Theorem}
\begin{proof}
By Lemma \ref{lem:hibi2}, $\bP_{\Delta(P)}$ has at worst terminal singularities. Then we can use the same formulas as in \cite[Proposition~8.6]{MR1463173}.
By substituting $d_j \Delta(P)$ for $\Delta_i$, we can derive the formulas~(\ref{eq:formula}). The explicit form of $\Delta_i$ is given by using a nef-partition in Section~\ref{subsub}.
\end{proof}

\subsection{Toric degenerations and conifold transitions}\label{sec:toric}
Here we summarize a key theorem by Gonciulea and Lakshmibai that gives a relationship between minuscule Schubert varieties and Hibi toric varieties.
\begin{Theorem}[{\cite[Theorem 7.34]{MR1417711}}]\label{thm:ToricDegen}
The minuscule Schubert variety $X(w)$ degenerates to the Hibi toric variety $\bP_{\Delta(P_w)}$, where $P_w$ is the minuscule poset for $X(w)$. More precisely, there exists a flat family $\scV \rightarrow \bC$ such that $\scV_t \simeq X(w)$ for all $t \in \bC^*$ and $\scV_0 \simeq \bP_{\Delta(P_w)}$.
\end{Theorem}

\begin{Corollary}\label{cor:terminal} A Gorenstein minuscule Schubert variety has at worst terminal singularities.
\end{Corollary}
\begin{proof}
By Lemma \ref{lem:hibi2}, a Gorenstein Hibi toric variety has at worst terminal (and hence canonical) singularities. It is known that a~deformation of a~terminal singularity is again terminal \cite{Nakayama1998} as with the case of canonical singularities~\cite{MR1631527}. Hence the assertion follows from Theorem~\ref{thm:ToricDegen}.
\end{proof}

A \textit{conifold transition} of a smooth Calabi--Yau 3-fold $X$ is the composite operation of a f\/lat degeneration of $X$ to $X_0$ with f\/initely many
nodes and a small resolution $Y\rightarrow X_0$. For a~general complete intersection Calabi--Yau 3-fold $X_0$ in a Gorenstein Hibi toric variety $\bP_{\Delta(P)}$, $X_0$ has at worst f\/initely many nodes by a Bertini type theorem for toroidal singularities. In fact, 3-dimensional Gorenstein terminal toric singularities are at worst nodal since they are presented by 3-dimensional cones over a unit triangle or a unit square up to unimodular transformations.
Thus, if a smooth Calabi--Yau 3-fold $X$ degenerates into $X_0 \subset \bP_{\Delta(P)}$, we obtain a conifold transition $Y$ of $X$, which is a strict transform of $X_0$ by a MPCP resolution $\hat \bP_{\Delta(P)}\rightarrow \bP_\Delta(P)$. Thus, Corollary~\ref{cor:coni} is a direct consequence of Theorem~\ref{thm:ToricDegen}.

\begin{Corollary}\label{cor:coni}
Let $X$ be a smooth complete intersection Calabi--Yau $3$-fold in a Gorenstein minuscule Schubert variety $X(w)$. There exists a conifold transition $Y$ of $X$ such that~$X$ degenera\-tes to $X_0 \subset \bP_\Delta(P_w)$ and~$Y$ is a strict transform of~$X_0$ by a MPCP resolution $\hat \bP_{\Delta(P_w)}\rightarrow \bP_\Delta(P_w)$, where $P_w$ is the minuscule poset for $X(w)$.
\end{Corollary}

\section{List of complete intersection Calabi--Yau 3-folds}\label{sec:list}
In this section, we study smooth complete intersection Calabi--Yau 3-folds in minuscule Schubert varieties. We show that there are three new deformation equivalent classes of such Calabi--Yau 3-folds. Among them, we f\/ind only one new example of Calabi--Yau 3-folds of Picard number one, that is, the complete intersection of nine hyperplanes in a locally factorial Schubert variety~${\boldsymbol{\Sigma}}$ def\/ined in Def\/inition~\ref{dfn:sigma}.

First we f\/ix some basic terminologies to describe our result. Let $X(w)$ be a minuscule Schubert variety. We call a subvariety $X \subset X(w)$ a~\textit{complete intersection} if it is the common zero locus of $r=\mathrm{codim} X$ global sections of invertible sheaves on~$X(w)$. We denote by $X=X(w)(d_1, \dots , d_r)$ the complete intersection variety of general $r$ sections of degree $d_1, \dots , d_r$ with respect to $\scO_{G/Q}(1)|_{X(w)}$,
the ample generator of $\Pic X(w)\simeq \bZ$. A \textit{Calabi--Yau} variety $X$ is a~normal projective variety with at worst Gorenstein canonical singularities and with trivial canonical bundle $K_X\simeq 0$ such that $H^i(X,\scO_X) = 0$ for all $0< i < \dim X$. Two smooth varieties~$X_1$ and~$X_2$ are called \textit{deformation equivalent} if there exists a~smooth family $\scX \rightarrow U$ over a~connected open base $U \subset \bC$ such that $\scX_{t_1}\simeq X_1$ and $\scX_{t_2}\simeq X_2$ for some $t_1, t_2 \in U$.

In the following proposition, we have all possible smooth complete intersection Calabi--Yau 3-folds in minuscule Schubert varieties:
\begin{Proposition}\label{prop:list}
A smooth complete intersection Calabi--Yau $3$-fold in a minuscule Schubert variety is one of those listed in Table~{\rm \ref{tab:smooth}} up to deformation equivalence. The Picard numbers of these Calabi--Yau $3$-folds are as listed on the last line.
\begin{table}[h]\centering
 \includegraphics{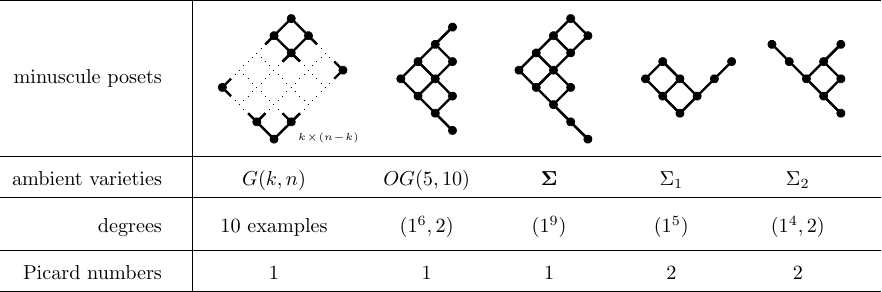}
\caption{Calabi--Yau 3-folds in minuscule Schubert varieties.}\label{tab:smooth}
\end{table}
In Table~{\rm \ref{tab:smooth}}, $10$ known examples in Grassmannians of type~A include five in projective spaces;
\[ \bP^4(5),\ \bP^5(2,4),\ \bP^5\big(3^2\big),\ \bP^6\big(2^2,3\big) \text{ and } \bP^7\big(2^4\big), \]
and five in others, whose mirror symmetry was discussed in~{\rm \cite{MR1619529}};
\[ G(2,5)\big(1^2,3\big),\ G(2,5)\big(1,2^2\big),\ G(2,6)\big(1^4,2\big), \ G(3,6)\big(1^6\big) \text{ and } G(2,7)\big(1^7\big).
\]
The Schubert varieties $\Sigma_1\subset G(3,7)$ and $\Sigma_2\subset OG(6,12)$ are defined by their minuscule posets in the first line.
\end{Proposition}

Before giving a proof of Proposition~\ref{prop:list}, we show a lemma on the Hibi toric varieties to handle the complexity of the redundant appearance of deformation equivalent varieties. A special case of this lemma was already used in~\cite{MR1619529} to prove the mirror duality of Hodge numbers for a~Calabi--Yau 3-fold $G(2,5)(1^2,3)$ and its mirror by reducing the argument to that of a hypersurface in a~Hibi toric variety.

\begin{Lemma}\label{lem:cone}Let $P$ be a finite poset and $P^*=P\cup \{ \hat{1} \}$, $P_*=P\cup \{ \hat{0} \}$ the posets	where $u \prec \hat{1}$, $\hat{0} \prec u$ for all $u \in P$, respectively. The Hibi toric varieties $\bP_{\Delta(P^*)}$ and $\bP_{\Delta(P_*)}$ are the projective cones over $\bP_{\Delta(P)}$ in $\Proj \bC[J(P^*)]$ and $\Proj \bC[J(P_*)]$, respectively.
\end{Lemma}
\begin{proof}	Let $1=1_{J(P^*)}$ be the maximal element in $J(P^*)$.	This $1$ is join irreducible in $J(P^*)$ because $J(P)$ has the unique maximal element. 	Therefore the variable $p_1$ in the homogeneous coordinate ring $A_{J(P^*)}$ of $\bP_{\Delta(P^*)}$	does not appear in any relation $p_\tau p_\phi - p_{\tau\wedge\phi} p_{\tau\vee\phi}$ $( \tau \not\sim \phi)$. Thus the Hibi toric variety $\bP_{\Delta(P^*)}$ is the projective cone over $\bP_{\Delta(P)}$.	The same argument is also valid for the minimal element $0=0_{J(P_*)}$ in $J(P_*)$ by considering the order dual to~$P$.
\end{proof}

\begin{proof}[Proof of Proposition \ref{prop:list}] We may assume that the ambient minuscule Schubert variety is Gorenstein. In fact, from the adjunction formula and the Grothendieck--Lefschetz theorem for divisor class groups of normal projective varieties~\cite{MR2219849}, we have an explicit formula of the canonical divisor as a~Cartier divisor, $K_{X(w)}=-D_1-\cdots -D_r$ where $D_j\subset X(w)$ is a very ample Cartier divisor of degree $d_j$ and $X(w)(d_1, \dots , d_r)$ is a~general Calabi--Yau complete intersection.

Let $X(w)$, $X(w')$ be minuscule Schubert varieties and $P_{w}$, $P_{w'}$ the corresponding minuscule posets, respectively. Assume that $P_{w'}$ coincides with a $d$-times iterated extension $(P_w)^{*\cdots *}_{*\cdots *}$ of~$P_w$. From Lemma~\ref{lem:cone}, it holds that $\bP_{\Delta(P_w)}$ is isomorphic to a complete intersection of $d$ general hyperplanes in~$\bP_{\Delta(P_{w'})}$. By Theorem~\ref{thm:ToricDegen}, there exist toric degenerations of~$X(w)$ and~$X(w')$ to the Hibi toric varieties $\bP_{\Delta(P_w)}$ and $\bP_{\Delta(P_{w'})}$, respectively. This means that general complete intersection Calabi--Yau 3-folds $X=X(w)(d_1, \dots, d_r)$ and $X'=X(w')(1^d, d_1, \dots, d_r)$ can be connected by f\/lat deformations through a complete intersection $X_0=\bP_{\Delta(P_w)}(d_1,\dots, d_r)$. Since~$X_0$ has at worst terminal singularities, the Kuranishi space is smooth by \cite[Theorem~A]{MR1286924} and the degenerating locus have a positive complex codimension. Therefore $X$ and $X'$ are connected by smooth deformation. Thus we can eliminate such redundancy arisen from iterated extensions of minuscule posets.

A Gorenstein minuscule Schubert variety $X(w)$ with minuscule poset $P=P_w$ is a $\abs{P}$-dimensional Fano variety of index $h_P$ as we stated in Proposition~\ref{prop:schubertsing}(1). The condition that general complete intersections in $X(w)$ are Calabi--Yau 3-folds gives a strong combinatorial restriction for the poset $P$ as follows,
\begin{gather*} h_P - 1 \le |P| \le h_P +3.\end{gather*}
On the other hand, there is a complete list of the minuscule posets in~\cite{MR2360316}. Hence we can make a list of the complete intersection Calabi--Yau 3-folds by counting such posets.

In fact, there are inclusions of Grassmannians of type~A and type~D as Schubert varieties,
\begin{gather*}
G(k, n-1) \subset G(k,n),\\
G(k-1,n-1)\subset G(k,n),\\
OG(m-1,2m-2) \subset OG(m,2m)
\end{gather*}
for each $k$, $n$ and $m$. The inclusions are compatible with the action of the Borel subgroups. Then it is enough to check the examples with $k \le 5$, $n-k\le 5$ and $m \le 7$. Otherwise, any Schubert variety that satisf\/ies the condition $\abs P \le h_P +3$ is included in smaller Grassmannians as their Schubert variety. The examples in Grassmannians of other types are at most f\/inite.

By forgetting the dif\/ference of projective cones by Lemma~\ref{lem:cone}, one can see that almost all the complete intersections obtained in this way are nothing but complete intersections in $G(k,n)$, $OG(5,10)$ or a~direct product of projective spaces. Other newly obtained examples correspond to the following f\/ive minuscule posets. The (essential) holes are circled in the following diagrams:
\begin{center}
\includegraphics{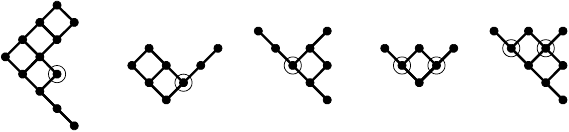}
\end{center}
In these Schubert varieties $X(w)$, the complete intersection 3-folds $X$ with trivial canonical bundles turn out to be Calabi--Yau varieties by some computation using the vanishing theorems for $X(w)$ given in Theorem \ref{thm:vanishing}.

One can determine whether these Calabi--Yau 3-folds are smooth or not by looking at the codimension of the singular locus of $X(w)$ using Proposition~\ref{prop:schubertsing}(2). For example, a general linear section $X={\boldsymbol{\Sigma}}(1^9)$ is smooth since the singular locus of ${\boldsymbol{\Sigma}}$ has codimension~7 as we saw in Proposition~\ref{prop:bSigma}. There are two more smooth examples.

For smooth complete intersections, the Picard number equals the number of peaks of minuscule posets since the number is the rank of the divisor class group $\mathrm{Cl}(X(w))$ again by the Grothendieck--Lefschetz theorem for divisor class groups~\cite{MR2219849}. This completes the proof.
\end{proof}

\section{Topological invariants}\label{sec:invariants}
In this section, we describe our calculation of topological invariants for smooth complete intersection Calabi--Yau 3-folds in minuscule Schubert varieties by taking $X={\boldsymbol{\Sigma}}(1^9)$ as an example. The topological invariants mean the degree $\deg(X) = \int_X H^3$, the linear form associated with the second Chern class $c_2(X)\cdot H= \int_X c_2(X)\cup H$ and the Euler number $\chi(X)=\int_X c_3(X)$, where $H$ is the ample generator of $\Pic(X) \simeq \bZ$.
\begin{Proposition}\label{prop:invariants}
The topological invariants of $X={\boldsymbol{\Sigma}}(1^9)$ are
\begin{gather*} \deg (X)=33,\qquad c_2(X)\cdot H=78,\qquad \chi(X)=-102.\end{gather*}
\end{Proposition}

\begin{proof} The degree of $X$ coincides with that of the minuscule Schubert variety ${\boldsymbol{\Sigma}} \subset \mathbb{OP}^2$ since the ample generator $\scO_{{\boldsymbol{\Sigma}}}(1)$ of $\Pic {\boldsymbol{\Sigma}}$ is the restriction of $\scO_{\mathbb{OP}^2}(1)$ and~$X$ is a linear section. We obtain $\deg ({\boldsymbol{\Sigma}})=33$ by using the Chevalley formula of $\bO\bP^2$ as we already saw in Remark~\ref{rem:minuscule}.

The Schubert variety $V^0:={\boldsymbol{\Sigma}}$ and its general complete intersections $V^j:={\boldsymbol{\Sigma}}(1^j)$ have at worst rational singularities (in fact at worst terminal singularities by Corollary \ref{cor:terminal}). Hence the Kawamata--Viehweg vanishing theorem gives
	\begin{gather*}
					H^i\big(V^j,\omega_{V^j} \otimes \scO_{V^j}(k)\big)=H^i\big(V^j, \scO_{V^j}(k+j-h_\bfP)\big)=0 \text{\quad for all } i > 0 \text{ and } k > 0. \end{gather*}
 Together with the long cohomology exact sequences from
	\begin{gather*}
					0 \rightarrow \scO_{V^j}(k) \rightarrow \mathcal{O}_{V^j}(k+1) \rightarrow \mathcal{O}_{V^{j+1}}(k+1) \rightarrow 0, \end{gather*}
	we calculate the holomorphic Euler number of $X=V^{h_\bfP}$ as
	\begin{gather*} \chi(X,\mathcal{O}_X(1))=\dim H^0(X,\mathcal{O}_X(1))=\dim H^0({\boldsymbol{\Sigma}},\mathcal{O}_{{\boldsymbol{\Sigma}}}(1))-h_\bfP= \abs{J(P)}-9=12. \end{gather*}
 On the other hand, it holds that
	\begin{gather*}
	\chi(X,\mathcal{O}_X(1))=\frac{1}{6}\deg (X)+\frac{1}{12}c_2(X)\cdot H
	\end{gather*}
from the Hirzebruch--Riemann--Roch theorem of the smooth Calabi--Yau 3-fold $X$. Thus we obtain $c_2(X)\cdot H=78$.

For the topological Euler number $\chi(X)$, we use the toric degeneration of ${\boldsymbol{\Sigma}}$ to the Hibi toric variety~$\bP_{\Delta(\bfP)}$ given in Theorem~\ref{thm:ToricDegen}. 	As we saw in Corollary~\ref{cor:coni}, we have a conifold transition~$Y$ of~$X$, which is a smooth Calabi--Yau complete intersection in a MPCP resolution~$\widehat{\bP}_{\Delta(\bfP)}$ of~$\bP_{\Delta(\bfP)}$. By Theorem~\ref{thm:hodge}, the Hodge numbers of~$Y$ can be calculated as $h^{1,1}(Y)=5$ and
 \begin{gather*}
h^{2,1}(Y)=9 \left(\abs{J(\bfP)}-9 \right)-\sum_{e\in E}(l^*(9\theta_e)-9l^*(8\theta_e))-\abs \bfP\\
\hphantom{h^{2,1}(Y)}{} =96-\sum_{e\in E}(l^*(9\theta_e)-9l^*(8\theta_e)).
 \end{gather*}
To count the number of interior integral points in each facet, we use Proposition \ref{prop:facestr} which states that a~face of the order polytope is also the order polytope of some poset~$\bfP'$. For each facet $\theta_e$, the corresponding poset $\bfP'$ (or $\hat{\bfP}'$) is easily obtained by replacing the inequality $x_{s(e)} \ge x_{t(e)}$ by the equality $x_{s(e)} = x_{t(e)}$ and by considering the induced partial order. The Hasse diagrams of resulting posets $\bfP'$ are shown in the following table, where the numbering of edges is chosen from the upper left of the Hasse diagram of $\hat{\bfP}$ like in Fig.~\ref{fig:hasse}.
\begin{table}[h]\centering
 \includegraphics{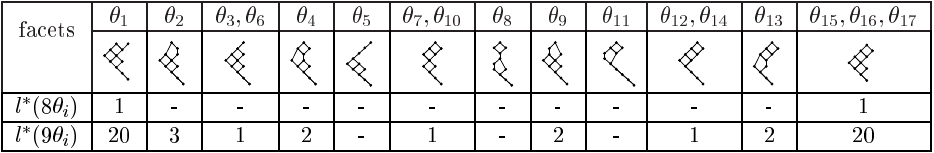}
\end{table}
As shown in this table, some posets~$\bfP'$ are pure and others are not. For a pure poset $\bfP'$, the face $\theta_e \simeq \Delta(\bfP')$ def\/ines Gorenstein Hibi toric variety~$\bP_{\Delta(\bfP')}$. Then we know $l^*(h_{\bfP'}\theta_e)=1$ and $l^*((h_{\bfP'}+1)\theta_e)=\abs{J(\bfP')}$. When $\bfP'$ is not pure, we can also easily obtain the number $l^*(k \theta_e)$ by counting the points satisfying the inequalities of the polytope $k\theta_e \simeq k \Delta(\bfP')$	strictly. For example, $9\theta_2$ contains three internal integral points corresponding to
\begin{table}[h]\centering
 \includegraphics{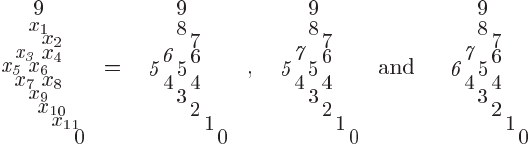}
\end{table}

\noindent
where only the integers written in italic can change under the inequalities. From the table, we obtain $h^{2,1}(Y)=37$, hence $\chi(Y)= 2(h^{1,1}(Y)-h^{2,1}(Y))= -64$.

Let us recall that the conifold transition is 	a surgery of Calabi--Yau 3-folds replacing f\/inite vanishing $S^3$ by the same number of exceptional $\bP^1\simeq S^2$ as in~\cite{MR690465}. From the inclusion-exclusion principle of the Euler numbers, $\chi(X)$ and $\chi(Y)$ are related with each other as
\begin{gather*} \chi(X)=\chi(Y)-2p,\end{gather*}
where $p$ is the number of nodes on $X_0$. Then we need to know the total degree of codimension three singular locus of~$\bP_{\Delta(\bfP)}$, which equals the number $p/\prod_j d_j=p$ in our case. From Theorem~\ref{thm:sing}, an irreducible component of singular locus of the Hibi toric variety~$\bP_{\Delta(\bfP)}$ corresponds to a minimal convex cycle in~$\hat \bfP$. There are four such cycles (or boxes) $b_1, \dots , b_4$ and all of them def\/ine the codimension three faces in~$\Delta(\bfP)$	as in Proposition~\ref{prop:facestr}. Again we can compute the corresponding index poset~$\bfP'$ of them by the method used above. The resulting posets $\bfP'$ are summarized as follows
\begin{table}[h]\centering
 \includegraphics{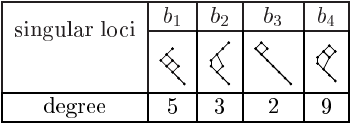}
\end{table}

From Proposition \ref{prop:deg}, we can compute the degree of each irreducible component of singular locus by counting the maximal chains in $J(\bfP')$. Then we obtain that total degree $p$, that is, the number of nodes on $X_0$ is $19$. We conclude $\chi(X)=-102$.
\end{proof}

\begin{Remark}\quad
 \begin{enumerate}\itemsep=0pt
\item[(1)] The existence of the Calabi--Yau 3-fold with these topological invariants was previously conjectured by \cite{MR2282974} from the monodromy calculations of Calabi--Yau dif\/ferential equations. We also perform a similar calculation in the next section.
\item[(2)] It may be possible to calculate the Euler number $\chi(X)$ in another way, by computing the Chern--Mather class of the Schubert variety~${\boldsymbol{\Sigma}}$. For the Grassmann Schubert varieties, this is done by~\cite{MR2628830} using Zelevinsky's $IH$-small resolution. In our case, however, it is known that ${\boldsymbol{\Sigma}}$ does not admit any $IH$-small resolution~\cite{MR2360316}.
 \end{enumerate}
\end{Remark}

\section{Mirror symmetry}\label{sec:mirror}
In this section, we study the mirror symmetry for Calabi--Yau complete intersections of Picard number one in minuscule Schubert varieties. In Section~\ref{ss:mirror}, we focus on the specif\/ic example $X={\boldsymbol{\Sigma}}(1^9)$ and carry out calculations related to mirror symmetry.

\subsection{Mirror constructions}
Let $X(w)$ be a Gorenstein minuscule Schubert variety and $X$ a Calabi--Yau complete intersection of Picard number one. Our purpose here is to describe a conjectural mirror family for~$X$ based on the general method proposed by~\cite{MR2112571, MR1619529} and compute the fundamental period of the family.

\subsubsection{Batyrev--Borisov construction} \label{subsub}
First we recall the Batyrev--Borisov mirror construction \cite{MR1269718,1993alg.geom.10001B} and apply it to Calabi--Yau complete intersections in Gorenstein Hibi toric varieties.

Let $\Delta$ and $\Delta^*$ be a polar dual pair of ref\/lexive polytopes. A \textit{nef-partition} is a~Minkowski sum decompositon $\Delta = \Delta_1 + \dots + \Delta_r$ into lattice polytopes $\Delta_j$ $(j=1,\dots , r)$ associated with nef line bundles. It is shown by~\cite{1993alg.geom.10001B} that a~nef-partition has a \textit{dual nef-partition} $\nabla=\nabla_1+ \cdots + \nabla_r$ in the dual vector space, def\/ined by
$\langle \Delta_i, \nabla_j \rangle \geq -\delta_{ij}$ for all $1\le i, j \le r$. A~general af\/f\/ine complete intersection
\begin{gather*}
X^\circ_{\Delta_1,\dots, \Delta_r}:= \bigg\{\sum_{m\in \Delta_j} a_m t^m =0 \, \Big| \, j=1,\dots,r \bigg\} \subset (\bC^*)^n
\end{gather*}
can be compactif\/ied to Calabi--Yau varieties $X_{\Delta_1,\dots, \Delta_r}$ and $\hat X_{\Delta_1,\dots, \Delta_r}$ in $\bP_\Delta$ and a MPCP resolution $\hat \bP_\Delta$ of $\bP_\Delta$, respectively. Here $X_{\Delta_1,\dots, \Delta_r}$ is nothing but a general complete intersection of corresponding nef line bundles on~$\bP_\Delta$. Thus one obtain the Batyrev--Borisov mirror pair of Calabi--Yau varieties, i.e., $\hat X_{\Delta_1,\dots, \Delta_r}$ and $\hat X_{\nabla_1,\dots, \nabla_r}$, based on the duality of nef-partitions.

Now we apply the Batyrev--Borisov mirror construction to Calabi--Yau complete intersections in Gorenstein Hibi toric varieties. Let $P$ be a f\/inite pure poset. Recall that the order polytope $\Delta(P) \subset M_\bR$ is associated with the hyperplane class on the Gorenstein Fano Hibi toric va\-rie\-ty~$\bP_{\Delta(P)}$ of index~$h_P$, and $\Delta= \sum\limits_{u\in P}h(u)\chi_u- h_P \Delta(P)$ is a corresponding ref\/lexive polytope to $-\scK_{\bP_{\Delta(P)}}= \scO(h_P)$.

There is a convenient description of the polar dual polytope $\Delta^*$ (cf.\ \cite{MR1756568,MR2770546}), where the set of edges~$E$ in the Hasse diagram of $\hat{P}$ plays an important role again. We regard the abelian groups $\bZ \hat{P}=N\oplus \bZ \{\hat{0}, \hat{1}\}$ and~$\bZ E$ as the groups of $0$-chains and $1$-chains of the natural chain complex associated with the Hasse diagram of $\hat{P}$. Here the boundary map in the chain complex is
\begin{gather*} \partial\colon \ \bZ E \longrightarrow \bZ \hat{P}, 	\qquad e \mapsto t(e)- s(e). \end{gather*}
We also consider the projection $\mathrm{pr}_1\colon \bZ \hat{P} \rightarrow N$ and the composed map
\begin{gather*}
\delta:=\mathrm{pr}_1 \circ \partial\colon \ \bZ E \longrightarrow N.\end{gather*}
One can see that the dual polytope $\Delta^*$ coincides with the convex hull of the image $\delta(E)\subset N_\bR$ directly from the def\/initions. The linear map $\delta$ gives a bijection between $E$ and the set of vertices in $\Delta^*$.

Let $X_0 \subset \bP_{\Delta(P)}$ be a general Calabi--Yau complete intersection of degree $(d_1, \dots , d_r)$ with respect to~$\scO(1)$. That is, $d_1, \dots, d_r$ satisf\/ies $\sum\limits_{j=1}^r d_j=h_P$. Let us choose a \textit{nef-partition} of~$\Delta$ in the following specif\/ic way. Def\/ine subsets~$E_j$ of edges in~$E$ as
\begin{gather*}
E_j= \bigcup_{k=d_1+\cdots+d_{j-1}+1}^{d_1+\cdots+d_j}E^k,
\end{gather*}
where $E^k:= \{e\in E \,|\, h(s(e))=k \}$. Then it holds $\scO(D_{E_j})=\scO(d_j)$ and we obtain a nef-partition $\Delta = \Delta_1 + \dots + \Delta_r$ with $\nabla_j=\Conv(\{0\}, \delta(E_j))$ and $\left< \Delta_i, \nabla_j \right> \geq -\delta_{ij}$. Let us set $\nabla=\nabla_1+ \cdots + \nabla_r \subset N_\bR$, which gives another ref\/lexive polytope~$\nabla$ and the dual nef-partition.

\begin{Proposition}[{cf.\ \cite{1993alg.geom.10001B}}] Let $P$ be a finite pure poset.	Denote by $X_0$	a complete intersection Calabi--Yau variety	of degree $(d_1, \dots , d_r)$ in a Gorenstein Hibi toric variety $\bP_{\Delta(P)}$,	and by \smash{$Y=\hat{X_0}$} the strict transform of $X_0$ in a MPCP-resolution $\hat{\bP}_{\Delta(P)}$ of $\bP_{\Delta(P)}$. A~Batyrev--Borisov mirror of~$Y$ is birational to the variety $X_{\nabla_1,\dots, \nabla_r}^\circ$ defined by the following equations in $(\bC^*)^{|P|}$:
\begin{gather*}
\tilde{f}_j =1-\bigg(\sum_{e\in E_j}a_e t^{\delta(e)}\bigg)=0 \quad \text{for all $1 \le j \le r$},\end{gather*}
where $a_e\in \bC$ are parameters. 	Further, the mirror Calabi--Yau variety	$Y^*$ of $Y$ is also obtained as the closure of the above set in	a~MPCP-resolution $\hat{\bP}_{\nabla}$ of $\bP_\nabla$.
\end{Proposition}
The mirror $Y^* \subset \hat{\bP}_{\nabla}$ actually has the expected stringy (or string-theoretic) Hodge numbers as proved in general by \cite{MR1408560, MR1463173}.

\subsubsection{Construction via conifold transition}\label{conifold}

\begin{figure}[h]\centering
 \includegraphics{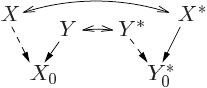}
 \caption{Mirror symmetry and conifold transitions.}\label{Fig:conifold}
\end{figure}

A conjectural mirror construction via conifold transitions of smooth Calabi--Yau 3-folds is proposed by \cite{MR1619529}. Let $X$ be a smooth Calabi--Yau 3-fold and $Y$ another Cabi--Yau 3-fold obtained by a conifold transition through the Calabi--Yau 3-fold $X_0$ with f\/initely many nodes. Morrison conjectured that the mirror Calabi--Yau 3-folds $Y^*$ and $X^*$ are again related in the same way~\cite{MR1673108}, see Fig.~\ref{Fig:conifold}. In the diagram, dashed and solid arrows represent f\/lat degenerations and small contraction morphisms respectively.
By an argument in \cite{MR2112571} on generalized monomial-divisor correspondence, there is a natural specialization~$Y^*_0$ of the family of~$Y^*$ to get the mirror of~$X$. That is, the specialized parameter $(a_e)_{e\in E}$ should be $\Sigma(\Delta^*)$-\textit{admissible}, i.e., there exists a~$\Sigma(\Delta^*)$-piecewise linear function $\phi\colon N_\bR \rightarrow \bR$ corresponding to a~Cartier divisor on~$X$ such that $\phi\circ \delta(e)=\log |a_e|$. In our cases, $\Pic \bP_{\Delta(P)} \simeq \Pic X\simeq \bZ$ holds. Then we can simply specialize the family to be diagonal, i.e., setting all the coef\/f\/icients~$a_e$ to be the same parameter~$a$. Now we present the conjecture of~\cite{MR1619529} in a~particular case with our terminology and
combining it with the argument on the number of nodes in the proof of Proposition~\ref{prop:invariants}.
\begin{Conjecture}[{cf.\ \cite[Conjecture 6.1.2]{MR1619529}}]\label{conj:mirror}Let $P$ be a finite pure poset and $X_0$ a complete intersection Calabi--Yau $3$-fold in $\bP_{\Delta(P)}$. Assume that a flat degeneration from a smooth Calabi--Yau $3$-fold $X$ of Picard number one to~$X_0$. We define a one-parameter family of affine complete intersections in~$(\bC^*)^{\abs P}$ by the following equations:
\begin{gather}	\label{eq:mirror}
	f_j =1-a\bigg(\sum_{e\in E_j} t^{\delta(e)}\bigg)=0	\quad \text{for all $1 \le j \le r$}.
\end{gather}
The closure $Y_0^*$ of the above set in a MPCP-resolution $\hat{\bP}_{\nabla}$ has $p$ nodes, satisfying $(|E|-|P|-p-1)$ relations for general $a\in \bC$, where $p$ is a number of nodes on a Calabi--Yau $3$-fold $X_0\subset \bP_{\Delta(P)}$, namely, $p=(\prod_{j=1}^r d_j) \sum c_{\abs{P'}}$ where the sum is over all minimal convex cycles in~$P$ corresponding to the contractions $P \rightarrow P'$ with $\abs{P'}=\abs P -3$ and $c_{\abs{P'}}$ denotes the number of maximal chains in~$J(P')$. A~small resolution $X^* \rightarrow Y_0^*$ is a mirror manifold of~$X$ with the correct Hodge numbers, $h^{i,j}(X^*)=h^{3-j,i}(X)$.
\end{Conjecture}
\begin{Remark} In the case of smoothing of 3-dimensional Calabi--Yau hypersurfaces in Gorenstein Hibi toric varieties, we can see Conjecture~\ref{conj:mirror} holds by the same argument as in~\cite{MR1619529, MR2801412}, i.e., the polar duality of faces gives a one-to-one correspondence between singular $\bP^1 \subset \bP_{\Delta(P)}$ and torus invariant $\bP^1\times \bP^1 \subset \bP_{\Delta^*}$ which intersects non-transversally with the closure of the set $\{f_1=0 \}$, and the MPCP-resolution $\hat{\bP}_{\Delta^*}\rightarrow \bP_{\Delta^*}$ increases them by $h_P$ times.
\end{Remark}

\subsubsection{Fundamental period}
Let $X=X(w)(d_1, \dots, d_r)$ be a complete intersection Calabi--Yau variety of arbitrary dimension in minuscule Schubert variety~$X(w)$. Def\/ine a~one-parameter family of af\/f\/ine complete intersections by the same equations as~(\ref{eq:mirror}). We call this family a \textit{conjectural mirror family} for $X$ as well. Obviously, the coordinate transformation $t_u \rightarrow \zeta^{h_u} t_u$ gives a $\bZ_{h_P}$-symmetry $a\rightarrow \zeta a$ in this family, where $\zeta=e^{2\pi \sqrt{-1}/h_P}$. Therefore we should take $x:=a^{h_P}$ as a genuine moduli parameter.

The \textit{fundamental period} $\omega_0(x)$ of the conjectural mirror family is the normalized monodromy invariant period of the family around $x=0$. In the case that there is a smooth mirror $X^*$, $\omega_0(x)$ should coincide with an integration of the unique holomorphic $(|P|-r)$-form on $X^*$ over the monodromy invariant (real) torus cycle. For conjectural mirror family, we def\/ine and compute~$\omega_0(x)$ by the residue theorem as follows:
\begin{gather*}
\omega_0(x) := \frac{1}{(2\pi \sqrt{-1})^{|P|}}		\int_{|t_u|=1} \frac{1}{\prod\limits_{j=1}^r f_j} \prod_{i=1}^{|P|} \frac{\mathrm{d}t_i}{t_i}\\
\hphantom{\omega_0(x)}{}	= \sum_{m=0}^\infty a^{h_P m}\frac{1}{(2\pi \sqrt{-1})^{|P|}}
		\int_{|t_u|=1} \prod_{j=1}^r \bigg(\sum_{e\in E_j} t^{\delta(e)}\bigg)^{d_j m}\prod_{i=1}^{|P|} \frac{\mathrm{d}t_i}{t_i}\\
\hphantom{\omega_0(x)}{} = \sum_{m=0}^\infty x^m N(m) = \sum_{m=0}^\infty x^m	\frac{\prod\limits_{j=1}^r (d_j m)!}{m!^{h_P}} n(m),
\end{gather*}
where $J_j(m):=\{ (j,i) \,|\,1 \le i \le d_j m \}$, $J^k(m):=\{ (k,i) \,|\,1 \le i \le m \} \subset \bN^2$ and
\begin{gather}
	N(m):=\# \bigg\{ \phi\colon
		 \bigcup_j J_j(m) \rightarrow E\,\Big|\, \phi(J_j(m))\subset E_j,\	\sum_{s(e)=u} \#\phi^{-1}(e)=\sum_{t(e)=u} \# \phi^{-1}(e) \bigg\},\nonumber\\
	n(m):=\# \bigg\{ \phi\colon \bigcup_k J^k(m) \rightarrow E\, \Big|\, \phi(J^k(m))\subset E^k,\ \sum_{s(e)=u} \#\phi^{-1}(e)
		=\sum_{t(e)=u} \# \phi^{-1}(e) \bigg\}.				\label{eq:binom}
\end{gather}
The coef\/f\/icients $N(m)$ and $n(m)$ are counting the way to distribute the exponents to each monomial.

In the case that the Hasse diagram of $P$ (and hence $\hat{P}$) is a planar graph, we can simplify the formula further like~\cite{MR1619529}. This is originally formulated in the work of Bondal and Galkin \cite{Bondal-Galkin2010} for the Landau--Ginzburg mirror of minuscule $G/Q$. Note that a minuscule poset is always a~planar graph.

\begin{Definition}Let $P$ be a f\/inite pure poset whose Hasse diagram is a planar graph.
\begin{enumerate}\itemsep=0pt
\item[(1)] Def\/ine the dual graph $B$ of the Hasse diagram of $\hat{P}$ on a sphere~$S^2$ by putting both~$\hat{1}$ and~$\hat{0}$ on the same point. More presisely, the vertices of $B$ is the set of areas partitioned by edges of~$\hat P$, and the edges of $B$ corresponds to the edges~$E$ of~$\hat P$.
\item[(2)] Denote by $b_L, b_R \in B$ the farthest left and right areas respectively. Under this def\/inition of left and right, we def\/ine the orientation of an edge~$e$ of~$B$ as the direction from the left area $\bfl(e)$ to the right area $\bfr(e)$.
		\end{enumerate}
\end{Definition}

\begin{figure}[h]\centering
 \includegraphics{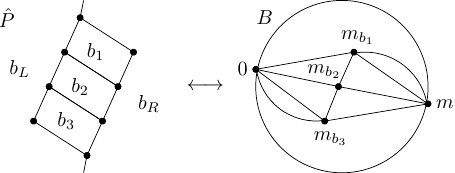}
\caption{An example of $\hat{P}$ and the dual graph $B$.}\label{Fig:B}
\end{figure}

An example of the Hasse diagram of $\hat P$ and its dual graph $B$ are drawn in Fig.~\ref{Fig:B}. Here we attach an indeterminate $m_b$ for each vertices $b\in B$ with $m_{b_L}=0$ and $m_{b_R}=m$, where $m$ is also an indeterminate.

\begin{Proposition} \label{prop:period}
Let $X(w)$ be a minuscule Schubert variety and $P=P_w$ its minuscule poset. Under the above notation, the fundamental period $\omega_0(x)$ for the conjectural mirror family of a~Calabi--Yau complete intersection $X(w)(d_1,\dots,d_r)$ is presented in the following:
	\begin{gather}\label{eq:binom2}
	\omega_0(x)=\sum_{m=0}^\infty \frac{\prod\limits_{j=1}^r (d_j m)!}{m!^{h_P}}
 \sum_{m_b\, (b \in B)} \prod_{e\in E} \begin{pmatrix} m_{\bfr(e)}\\m_{\bfl(e)} \end{pmatrix}x^m.
	\end{gather}
\end{Proposition}
\begin{proof}Let us compute $n(m)$ in (\ref{eq:binom}) for any poset $P$ whose Hasse diagram is a planar graph. Fix positive integers $m$ and $m_b\ (b\in B)$ where $m_{\bfl(e)} \le m_{\bfr(e)}$ for all $e \in E$ and $m_{b_L}= 0$ and $m_{b_R} =m$. Consider a map $\phi\colon \bigcup_k J^k(m) \rightarrow E$ satisfying the conditions $\phi(J^k(m))\subset E^k$ for all $1 \le k\le r$ and $\# \phi^{-1}(e) = m_{\bfr(e)}-m_{\bfl(e)}$. It also holds $\sum\limits_{s(e)=u} \#\phi^{-1}(e) =\sum\limits_{t(e)=u} \# \phi^{-1}(e)$ for all $u \in P$. It is clear that the number of such maps is $\prod\limits_{e\in E} \begin{pmatrix} m_{\bfr(e)}\\m_{\bfl(e)} \end{pmatrix}$. Conversely, we can determine such the set of positive integers $m_b$ by counting $\# \phi^{-1}(e)$ for each $e \in E$. This gives an explicit formula of~$n(m)$, and hence, the fundamental period $\omega_0(x)$ in~(\ref{eq:binom2}).
\end{proof}

\subsection[Mirror symmetry for ${\Sigma}(1^9)$]{Mirror symmetry for ${\boldsymbol{\Sigma}}(1^9)$} \label{ss:mirror}

\subsubsection{Picard--Fuchs equation} Now we focus on the Calabi--Yau 3-fold $X={\boldsymbol{\Sigma}}(1^9)$. The fundamental period of the mirror family
can be read from the following diagram
\begin{figure}[h]\centering
 \includegraphics{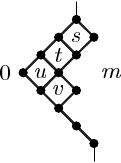}
\end{figure}

\noindent
The vertices of the dual graph $B$ corresponds to the separated areas. The fundamental period turns out to be
\begin{gather*}
		\omega_0(x)=\sum_{m, s, t, u, v}
		\begin{pmatrix} m \\ s \end{pmatrix}^2
		\begin{pmatrix} m \\ v \end{pmatrix}^2
 		\begin{pmatrix} m \\ t \end{pmatrix}
		\begin{pmatrix} s \\ t \end{pmatrix}
		\begin{pmatrix} t \\ u \end{pmatrix}
		\begin{pmatrix} v \\ u \end{pmatrix} x^m,
\end{gather*}
where $x=a^9$. With the aid of numerical method, we obtain the Picard--Fuchs equation for the conjectural mirror of~$X$.
\begin{Proposition}\label{prop:PF}
Let $\omega_0(x)$ be the above power series around $x=0$, which corresponds to the fundamental period for the conjectural mirror family of the Calabi--Yau $3$-fold $X={\boldsymbol{\Sigma}}(1^9)$. It satisfies the Picard--Fuchs equation $\mathcal{D}_x\omega_0(x)=0$ with $\theta_x=x \partial_x$ and
\begin{gather*}
 \mathcal{D}_x=121\theta_x^4 - 77x\big(130\theta_x^4 + 266\theta_x^3 + 210\theta_x^2 + 77\theta_x + 11\big)\\
 \hphantom{\mathcal{D}_x=}{} -x^2\big(32126\theta_x^4 + 89990\theta_x^3 + 103725\theta_x^2 + 55253\theta_x + 11198\big)\\
\hphantom{\mathcal{D}_x=}{} -x^3\big(28723\theta_x^4 + 74184\theta_x^3 + 63474\theta_x^2 + 20625\theta_x + 1716\big)\\
\hphantom{\mathcal{D}_x=}{} -7x^4\big(1135\theta_x^4 + 2336\theta_x^3 + 1881\theta_x^2 + 713\theta_x + 110\big) - 49x^5(\theta_x+1)^4.
\end{gather*}
\end{Proposition}

\begin{Remark}\quad\samepage
\begin{enumerate}\itemsep=0pt
	\item[(1)] The dif\/ferential operator $\scD_x$ in Proposition \ref{prop:PF} already appeared in the lists of \cite{2005math......7430A, CYdatabase} with conjectural topological invariants, obtained from the monodromy calculation for a~par\-ti\-cu\-lar basis of the solutions for $\scD_x$ and some assumptions by mirror symmetry \cite{MR2282974}. The existence of the smooth Calabi--Yau 3-fold $X={\boldsymbol{\Sigma}}(1^9)$ and our computation of topological invariants in Proposition \ref{prop:invariants} answer to their conjecture.
 \item[(2)] There exist two maximally uniportent monodromy points for $\scD_x$ in Proposition \ref{prop:PF}, not only $x = 0$ but $x = \infty$. Once we assume the existence of a suitable Calabi--Yau 3-fold $Z$ for $x=\infty$, one can compute conjectural topological invariants of $Z$ by the monodromy calculation based on the mirror symmetry conjecture \cite{MR1115626}. As with $X$ for $x=0$, this is already done by \cite{2005math......7430A, CYdatabase}.
\end{enumerate}
\end{Remark}
Thus we are led to the following conjecture based on the homological mirror symmetry similar to the examples of the Grassmannian--Pfaf\/f\/ian in~\cite{MR1775415} and the Reye congruence Calabi--Yau 3-fold in \cite{MR3166392}.
\begin{Conjecture}[{cf.\ \cite{2005math......7430A, CYdatabase}}] \label{conj:FM}
There exists a smooth Calabi--Yau $3$-fold $Z$ of Picard number one, whose derived category of coherent sheaves is equivalent to that of the Calabi--Yau $3$-fold $X={\boldsymbol{\Sigma}}(1^9)$. The topological invariants of $Z$ are
	\begin{gather*}
					\deg (Z)=21,\qquad c_2(Z)\cdot H=66,\qquad \chi(Z)=-102,\end{gather*}
 where $H$ is the ample generator of the Picard group $\Pic(Z)\simeq \bZ$.
\end{Conjecture}
The Calabi--Yau 3-fold $Z$ in Conjecture~\ref{conj:FM} cannot be birational to $X$ because $h^{1,1}=1$ and $\deg (Z) \ne \deg (X)$, so that it should be a non-trivial Fourier--Mukai partner of~$X$.

\subsubsection{BPS numbers}
As a further consistency check in Conjecture~\ref{conj:FM} and an application of the mirror construction, we carry out the computation of BPS numbers by using the methods proposed by~\cite{MR1240687,MR1301851, MR1115626} based on the mirror symmetry conjecture. The BPS numbers $n_g(d)$ are related with the Gromov--Witten invariants~$N_g(d)$ by some closed formula \cite{1998hep.th...12127G}, hence we obtain the prediction for Gromov--Witten invariants of~$X$ and~$Z$ from these computations.
We omit all the details here and only present the results in Appendix~\ref{sec:bps}. For the details, one can get many references in now. Here we have followed \cite{MR2481271}, where a very similar example to ours, the Grassmannian--Pfaf\/f\/ian Calabi--Yau 3-fold, has been analyzed.

\appendix

\section[Other examples in minuscule $G/Q$]{Other examples in minuscule $\boldsymbol{G/Q}$}\label{sec:others}
We remark that there are some other examples of Calabi--Yau complete intersections of Picard number one in minuscule homogeneous spaces~$G/Q$. In these cases one may be able to know more details because the structure of the small quantum cohomology ring $QH^*(G/Q)$ of $G/Q$ has been studied intensively by recent researches \cite{MR1454400, MR2421317,MR2027200} for instance.
In particular, the quantum Chevalley formula \cite[Proposition~4.1]{MR2421317} for minuscule homogeneous spaces $G/Q$ is given in terms of the minuscule poset~$P_w$. From this formula, we can calculate the quantum dif\/ferential system of minuscule $G/Q$ and reduce this f\/irst-order system to a~higher-order dif\/ferential system by some computer aided calculation (see \cite{MR1619529,MR2391365} for instance). For example, the orthogonal Grassmannian~$OG(5,10)$ and the Cayley plane $\bO\bP^2$ have the quantum dif\/ferential operators annihilating the component associated with the fundamental class as follows:
\begin{gather*}
\begin{split}&
\theta_q^{11}(\theta_q-1)^5 - q\theta_q^5(2\theta_q+1)\big(17\theta_q^2 + 17\theta_q + 5\big)+ q^2, \\
& \theta_q^{17}(\theta_q -1)^9 -3q \theta_q^9(2\theta_q +1)(3\theta_q^2+3\theta_q+1)(15\theta_q^2+15\theta_q +4)-3q^2(3\theta_q+2)(3\theta_q+4),\end{split}\end{gather*}
respectively, where $\theta_q=q \partial_q$. By the quantum hyperplane section theorem for homogeneous spaces~\cite{MR1719555}, we can determine the genus zero Gromov--Witten potential for each example of Calabi--Yau complete intersections in minuscule $G/Q$. In fact this procedure also works for general rational homogeneous spaces $G/Q$ by replacing the quantum Chevalley formula by the general one obtained in~\cite{MR2072765} though we do not go into the details here.

On the other hand, we already know the conjectural one-parameter mirror family. Hence comparing both sides, we can prove the classical genus zero `mirror theorem' \cite{MR1408320,MR1621573} for these examples combining Conjecture~\ref{conj:mirror}. In particular, we can do this for $OG(5,10)(1^6,2)$, one of the examples of smooth complete intersection Calabi--Yau 3-folds of Picard number one in minuscule Schubert varieties in Proposition~\ref{prop:list}. The Picard--Fuchs operator for the conjectural mirror family of $OG(5,10)(1^6,2)$ coincides with the conjecture in~\cite{MR2282974},
\begin{gather*}
				\theta_x^4 - 2x(2\theta_x+1)^2\big(17\theta_x^2 + 17\theta_x + 5\big)+ 4x^2(\theta_x+1)^2(2\theta_x+1)(2\theta_x+3), \end{gather*}
where $\theta_x=x \partial_x$.

Examples of Calabi--Yau $d$-folds with $d \ne 3$ seem to be of some interest. The $K3$ surface $OG(5,10)(1^8)$ is an example which has a nontrivial
Fourier--Mukai partner \cite{MR1714828}, whose mirror symmetry are investigated in \cite{MR2047679}. The Picard--Fuchs operator for the mirror family
of $OG(5,10)(1^8)$ is
\begin{gather*} \theta_x^3-x(2\theta_x+1)\big(17\theta_x^2+17\theta_x+5\big)+x^2(\theta_x+1)^3, \end{gather*}
where $\theta_x=x\partial_x$. It has exactly the same property as the Picard--Fuchs operators of conjectural mirror families for $G(2,7)(1^7)$ and ${\boldsymbol{\Sigma}}(1^9)$, with two maximally unipotent monodromy points. Moreover this operator is also well-known since appeared essentially
in A\'pery's proof of irrationality of~$\zeta(3)$, similar to the Picard--Fuchs operator of the mirror family for $G(2,5)(1^5)$ in the $\zeta(2)$ case.
From this viewpoint, it may be valuable to remark that the conjectural mirror family of Calabi--Yau 4-fold $\bO\bP^2(1^{12})$ also has the following 5th-order
Picard--Fuchs operator,
\begin{gather*} \theta_x^5-3x(2\theta_x+1)\big(3\theta_x^2+3\theta_x+1\big)\big(15\theta_x^2+15\theta_x+4\big)-3x^2(\theta_x+1)^3(3\theta_x+2)(3\theta_x+4),\end{gather*}
where $\theta_x=x\partial_x$. This operator is related to the value $\zeta(4)$ \cite{MR2282972, MR2068773}. The relationship between the quantum dif\/ferential system of rational homogeneous spaces and zeta values are also numerically conf\/irmed in~\cite{Galkin2008} with the def\/inition of the A\'pery characteristic class.

\section{BPS numbers}\label{sec:bps}
\begin{table}[H]\centering\footnotesize
 \begin{tabular}{c |l l l l l}
\hline
 $d$ & $g=0$ & $g=1$ &$g=2$ &$g=3$&$g=4$ \\
\hline
1&252&0&0&0&0\\
2&1854&0&0&0&0\\
3&27156&0&0&0&0\\
4&567063&0&0&0&0\\
5&14514039&4671&0&0&0\\
6&424256409&1029484&0&0&0\\
7&13599543618&112256550&5058&0&0\\
8&466563312360&9161698059&7759089&0&0\\
9&16861067232735&645270182913&2496748119&151479&0\\
10&634912711612848&41731465395267&438543955881&418482990& $-$3708\\
11&24717672325914858&2557583730349461&56118708041940&217285861284&33975180
\end{tabular}
\caption{BPS numbers $n_g^X(d)$ of $X$.}
\end{table}
\begin{table}[H]\centering\footnotesize
\begin{tabular}{c |l l l l l}
\hline
$d$ & $g=0$ & $g=1$ & $g=2$ & $g=3$&$g=4$\\
\hline
1&387&0&0&0&0\\
2&4671&0&0&0&0\\
3&124323&1&0&0&0\\
4&4782996&1854&0&0&0\\
5&226411803&606294&0&0&0\\
6&12249769449&117751416&27156&0&0\\
7&727224033330&17516315259&33487812&252&0\\
8&46217599569117&2252199216735&15885697536&7759089&0\\
9&3094575464496057&265984028638047&4690774243470&13680891072&1127008\\
10&215917815744645750&29788858876065588&1053460470463461&9429360817149&12259161360\\
\end{tabular}
\caption{BPS numbers $n_g^Z(d)$ of $Z$.}
\end{table}

\subsection*{Acknowledgements} The author would like to express his deep gratitude to his supervisor Professor Shinobu Hosono for valuable suggestions and warm encouragement. He greatly appreciates many helpful discussions with Daisuke Inoue, Atsushi Kanazawa and Fumihiko Sanda at the seminars we had in University of Tokyo. He would also like to thank Yoshinori Gongyo, Takehiko Yasuda, Atsushi Ito and Taro Sano for useful comments to improve the work. The author thanks the anonymous referees for providing a number of valuable comments and in particular for pointing out the oversight of the examples of Picard number two in Proposition~\ref{prop:list}. Part of this paper was written at Mathematisches Institute Universit\"at T\"ubingen during his stay from October~1 to December~25, 2012. He was supported in part by Institutional Program for Young Researcher Overseas Visits by JSPS for this stay. It is a pleasure to thank Professor Victor Batyrev for valuable comments and creating a nice environment for the author.

\pdfbookmark[1]{References}{ref}
\LastPageEnding

\end{document}